\documentclass[a4paper,12pt,reqno]{amsart}

\usepackage{amsmath}
\usepackage{amssymb}
\usepackage{mathrsfs}

\begingroup\newtheorem{theorem}{Theorem}[section]
\newtheorem{lemma}[theorem]{Lemma}

\endgroup
\theoremstyle{definition}
\begingroup\newtheorem{definition}[theorem]{Definition}
\newtheorem{remark}[theorem]{Remark}
\newtheorem{example}[theorem]{Example}\endgroup
\numberwithin{equation}{section}

\newcommand{\MM}{{\mathbb M}}
\newcommand{\NN}{{\mathbb N}}

\newcommand{\RR}{{\mathbb R}}

\newcommand{\esssup}{\mathop{\rm ess\,sup}}
\newcommand{\dv}{{\rm div}}
\newcommand{\LM}[2]{\hbox{\vrule width.4pt \vbox to#1pt{\vfill
\hrule width#2pt height.4pt}}}
\newcommand{\LLL}{{\mathchoice {\>\LM{7}{5}\>}{\>\LM{7}{5}\>}{\,\LM{5}{3.5}\,}{\,\LM{3.35}{2.5}\,}}}

\title[]
{Quasistatic crack growth for a cohesive zone model with prescribed crack path}

\author[Gianni Dal Maso]{Gianni Dal Maso}\address[Gianni Dal Maso]{SISSA, Via Beirut 4, 34014 Trieste, Italy}
\email[Gianni Dal Maso]{dalmaso@sissa.it}\author[Chiara Zanini]{Chiara Zanini}
\address[Chiara Zanini]{SISSA, Via Beirut 4, 34014 Trieste, Italy}\email[Chiara Zanini]{zaninic@sissa.it}

\begin{document}

\begin{abstract} In this paper we study the quasistatic crack growth for a cohesive zone model.
We assume that the crack path is prescribed and we study the time evolution of the crack
in the framework of the variational theory of rate-independent processes.
\end{abstract}

\maketitle



\begin{section}{INTRODUCTION}\label{introd}

In this paper we present a variational model for quasistatic crack growth in the presence of a cohesive 
force exerted between the lips of the crack.

The evolution of the crack is governed by an energy which is the sum of three terms: 
the bulk energy of the uncracked part, the energy dissipated in the fracture process, 
and the work of the external loads. 
The main mathematical difficulty is given by the fact that the fracture energy depends on 
the opening of the crack. For this reason we cannot apply directly the tools developed so far 
in the applications to fracture mechanics of the theory of free discontinuity problems (see \cite{F-M}, \cite{DM-T}, 
\cite{DM-T-2}, \cite{Ch}, \cite{Fra-Lar}, \cite{DM-F-T-04}, \cite{DM-F-T-03}).

To simplify the mathematical difficulties, we assume that the crack path is prescribed, and we focus only on the time evolution. 
This allows us to consider very general bulk and crack energies, which may include constraints on the crack opening, related 
to the infinitesimal noninterpenetration of matter. 
The evolution of the crack
is defined (see Definition~\ref{iqe} below) in the framework of Mielke's approach to a variational theory of rate-independent 
processes (see \cite{Mie-Beijing}, \cite{M-M-03}).  

We prove an existence result for the quasistatic evolution, by approximating the continuous-time problem by discrete-time problems, 
for which the evolution is defined by solving incremental minimum problems. 
The irreversibility of the crack process leads to introduce an auxiliary time-dependent function $t\mapsto \gamma(t)$ 
(see Section~\ref{setting} below), defined on the prescribed crack path, 
which takes into account the local history of the crack up to time $t$. 
The main mathematical difficulty in the proof is the compactness of the approximating functions 
$t\mapsto \gamma_k(t)$. This is solved by introducing a new notion of convergence of functions related to the problem, 
with good compactness and semicontinuity properties.
\end{section}

\begin{section}{SETTING}\label{setting}

The reference configuration is a bounded open set $\Omega$
of $\RR^n$ with Lipschitz boundary $\partial\Omega$, 
which can be written as the union 
of two disjoint Borel sets $\partial_0\Omega$ and $\partial_1\Omega$, with $\mathscr{H}^{n-1}(\partial_0\Omega)>0$ and 
$\partial_1\Omega$ relatively open. Here and henceforth  $\mathscr{H}^{n-1}$
denotes the $(n-1)$-dimensional Hausdorff measure. 
On $\partial_0\Omega$, the Dirichlet part of the boundary,
we will assign the boundary deformation, while on $\partial_1\Omega$, the Neumann part of the boundary,
we will prescribe surface forces. 

We assume that the cracks are contained in a compact $C^1$-orientable $(n-1)$-dimensional manifold $M\subset \Omega$
with boundary $\partial M$, such that $\Omega \smallsetminus M$ is connected. 
Therefore it is reasonable to take the deformation $u$ as a function in the space $W^{1,p}(\Omega\smallsetminus M;\RR^m)$,
so that the essential discontinuity points of $u$ are contained in $M$. Although the natural choice is $m=n$, there are no 
mathematical difficulties in considering an arbitrary $m\geq 1$. The case $m=1$ is used in the study of antiplane shears.
The number $p> 1$ depends on the bounds on the energy density considered below.

We take into account prescribed time-dependent boundary deformations $t \mapsto\psi(t)$, with $\psi(t)\in W^{1,p}(\Omega;\RR^m)$, 
in the sense that for each time $t\in[0,T]$ we consider only deformations $u\in W^{1,p}(\Omega\smallsetminus M;\RR^m)$ such that
\begin{equation*}
u=\psi(t) \quad \text{on $\partial_0\Omega$},
\end{equation*}
where the previous equality has to be considered in the sense of traces.
We assume also that, as a function of time, $t\mapsto\psi(t)$ is absolutely continuous from $[0,T]$ into $W^{1,p}(\Omega;\RR^m)$.

Thus the time derivative $t\mapsto \dot{\psi}(t)$ belongs to the space
$L^1([0,T];W^{1,p}(\Omega;\RR^m))$ and its spatial gradient $t\mapsto \nabla\dot{\psi}(t)$
belongs to the space $L^1([0,T];L^p(\Omega;\MM^{m \times n}))$.

We assume that the uncracked part of the body is hyperelastic
and that its bulk energy relative to the deformation $u\in W^{1,p}(\Omega\smallsetminus M;\RR^m)$
is of the form 
\begin{equation*}
\int_{\Omega \smallsetminus M} W(x,\nabla u )\, dx, 
\end{equation*}
where $W(x,\xi)$ is a given Carath\'eodory function 
$W\colon(\Omega\smallsetminus M) \times {\mathbb M}^{m{\times}n} \to \RR$ such that
\begin{itemize}
\item[$(W_1)$] $\xi \mapsto W(x,\xi)$ is quasiconvex and $C^1$ for every
$x \in \Omega \smallsetminus M;$
\item[$(W_2)$] there are two positive constants $a_0,a_1 $   
and two nonnegative functions $b_0,b_1\in L^1(\Omega\smallsetminus M)$ such that
\begin{equation}\label{Whp}
a_0 \left |\xi\right|^{p}- b_0(x)\leq W(x,\xi) \leq a_1
\left |\xi\right|^{p} + b_1(x),
\end{equation} 
for every $(x,\xi) \in (\Omega\smallsetminus M) \times {\mathbb M}^{m \times n}$.
\end{itemize}
Since $\xi \mapsto W(x,\xi)$ is rank-one convex on 
$\MM^{m \times n}$ for every $x \in \Omega \smallsetminus M$, we can 
deduce from $(\ref{Whp})$ an estimate for the partial gradient of $W$ with respect to $\xi$,
$\partial_{\xi} W\colon (\Omega\smallsetminus M)\times \MM^{m \times n} \to \MM^{m \times n}$. 
More precisely, there are a positive constant 
$a_2$ and a nonnegative function $b_2 \in L^1(\Omega \smallsetminus M)$
such that
\begin{equation}\label{deWhp}
\left|\partial_{\xi}W(x,\xi)\right| \leq a_2\left|\xi\right|^{p-1}+b_2(x),
\end{equation}
for every $(x, \xi)\in (\Omega \smallsetminus M)\times \MM^{m\times n}$. 

To shorten the notation we introduce the function 
$\mathcal{W}\colon L^p(\Omega\smallsetminus M; \MM^{m\times n})\to \RR$ defined by
\begin{equation*}
\mathcal{W}(\Psi):= \int_{\Omega\smallsetminus M} W(x,\Psi)\, dx,
\end{equation*}
for every $\Psi \in L^p(\Omega\smallsetminus M; \MM^{m\times n})$. By $(\ref{Whp})$ and $(\ref{deWhp})$
the functional $\mathcal{W}$ is of class $C^1$ on $L^p(\Omega\smallsetminus M; \MM^{m\times n})$ and its 
differential 
$\partial\mathcal{W}\colon L^p(\Omega\smallsetminus M; \MM^{m\times n})\to L^q(\Omega\smallsetminus M; \MM^{m\times n})$,
$p^{-1}+q^{-1}=1$, is given by
\begin{equation*}
\langle \partial\mathcal{W}(\Psi),\Phi\rangle = 
\int_{\Omega\smallsetminus M} \partial_{\xi} W(x,\Psi){\,:\,}\Phi\;dx,
\end{equation*}
for every $\Phi$, $\Psi\in L^p(\Omega\smallsetminus M;\MM^{m\times n})$, where 
$\langle \cdot,\cdot\rangle$ 
denotes the duality pairing between the spaces $L^q(\Omega\smallsetminus M;\MM^{m\times n})$
and $L^p(\Omega\smallsetminus M;\MM^{m\times n})$, and $\partial_{\xi} W(x,\Psi){\,:\,}\Phi$ denotes the scalar product
between the two matrices $\partial_{\xi} W(x,\Psi)$ and $\Phi$.

By the assumptions on $W$, the functions
$\mathcal{W}$ and $\partial\mathcal{W}$ satisfy the following properties:
there are two positive constants $\alpha_0,\alpha_1 $
and two nonnegative constants $\beta_0,\beta_1$ such that
\begin{equation}\label{whp}
\alpha_0 \left \|\Psi\right\|^{p}_p- \beta_0\leq \mathcal{W}(\Psi) \leq \alpha_1
\left \|\Psi\right\|^{p}_p + \beta_1,
\end{equation} 
for every $\Psi\in L^{p}(\Omega\smallsetminus M;\MM^{m\times n})$, and there is a positive constant $\alpha_2$ such that 
\begin{equation}\label{dewhp}
\langle\partial\mathcal{W}(\Psi),\Phi\rangle\leq \alpha_2(1+\|\Psi\|_p^{p-1})\|\Phi\|_p,
\end{equation}
for every $\Psi$, $\Phi\in L^p(\Omega\smallsetminus M;\MM^{m\times n})$.

For a fixed time $t \in [0,T]$, we assume that the external time-dependent 
loads  $\mathscr{L}(t)$ belong to $(W^{1,p}(\Omega\smallsetminus M;\RR^m))'$, 
the dual space of $W^{1,p}(\Omega\smallsetminus M;\RR^m)$. The duality product $\langle\mathscr{L}(t),u\rangle $ is 
interpreted as the work done by the loads on the deformation $u$.

Let us fix an orientation of $M$ and let $u^{\oplus}$ be the trace of $u$ on the positive side of $M$, 
and $u^{\ominus}$ be the trace of $u$ on the negative side of $M$.
The most general form of the work done by the external loads is given by
\begin{equation}\label{def_L}
\begin{split}
\langle\mathscr{L}(t),u\rangle = &\int_{\Omega \smallsetminus M} f(t)\, u\, dx 
+ \int_{\Omega \smallsetminus M} H(t){\,:\,}\nabla u\, dx \,+\\
& + \int_{\partial_1 \Omega} g(t)\, u\, d\mathscr{H}^{n-1}
+ \int_{M} (g^{\oplus}(t)\, u^{\oplus}+ g^{\ominus}(t)\, u^{\ominus})\, d\mathscr{H}^{n-1},
\end{split}
\end{equation}
where $f(t)\in L^q(\Omega\smallsetminus M;\RR^m)$, $H(t)\in L^q(\Omega\smallsetminus M;\MM^{m\times n})$,  
$g(t)\in L^q(\partial_1 \Omega;\RR^m)$, $g^{\oplus}(t)$ and $g^{\ominus}(t)\in L^q(M;\RR^m)$, 
with $p^{-1}+q^{-1}=1$.
Actually the representation theorem for $(W^{1,p}(\Omega\smallsetminus M;\RR^m))'$ shows that it is enough to use just 
the terms of the first line of (\ref{def_L}).
The terms in the second line have been added in order to write in an explicit way 
the contribution of the surface forces acting on the Neumann part of the boundary and on one or both sides of $M$. 

With these assumptions we do not exclude the possibility that $H(t)$ could be discontinuous on $M$. Moreover, observe
that if $f(t),H(t),g(t),g^{\oplus}(t)$ and $g^{\ominus}(t)$ are sufficiently regular, then 
\begin{equation*}
f(t) - \hbox{ div } H(t)
\end{equation*}
plays the role of the volume forces on $\Omega\smallsetminus M$,
\begin{equation*}
g(t) + H(t)\nu 
\end{equation*}
plays the role of the surface forces on $\partial_1 \Omega$, and
\begin{equation*}
g^{\oplus}(t) - H^{\oplus}(t)\nu  \qquad \hbox{and} \qquad
g^{\ominus}(t) + H^{\ominus}(t)\nu 
\end{equation*}
play the role of the surface forces acting on the positive (respectively negative) side of $M$, where $\nu$ 
is the outer unit normal to $\partial(\Omega\smallsetminus M)$. We observe that, by our positions, $\nu$ turns out to be 
the inner normal on the positive side of $M$; this is why in the last formula we take the minus sign in front of~$H^{\oplus}(t)\nu$.

We assume that, as a function of time, $t\mapsto \mathscr{L}(t)$ is absolutely continuous from $[0,T]$ 
into $(W^{1,p}(\Omega\smallsetminus M;\RR^m))'$. Thus the time derivative $t\mapsto \dot{\mathscr{L}}(t)$ 
belongs to the space $L^1([0,T];(W^{1,p}(\Omega\smallsetminus M;\RR^m))')$. 
If $\mathscr{L}(t)$ is represented by (\ref{def_L}), then the absolute continuity of $t\mapsto \mathscr{L}(t)$ 
follows from the absolute continuity of the functions $t\mapsto f(t)$, $t\mapsto H(t)$, $t\mapsto g(t)$, $t\mapsto g^\oplus(t)$, 
and $t\mapsto g^\ominus(t)$.

If the deformation $u$ has a nonzero jump $[u]= u^{\oplus}-u^{\ominus}$ on $M$, then the body has a crack on (part of) $M$. 
More precisely the crack is given by the set 
\begin{equation*}
\{x\in M: [u](x)\neq 0\}.
\end{equation*}

Let us consider now the work done to produce a crack.
If we neglect for a moment the problem of irreversibility, we may assume that this work can be written in the form
\begin{equation*}
 \int_M \varphi(x,[u])\,d\mathscr{H}^{n-1},
\end{equation*}
where $\varphi\colon M \times \RR^m \to [0,+\infty]$ 
satisfies the following properties
\begin{itemize}
\item[$(\varphi_1)$] $\varphi$ is a Borel function;
\item[$(\varphi_2)$] $\varphi(x,0)=0$ for $\mathscr{H}^{n-1}$-a.e.\ $x\in M$;
\item[$(\varphi_3)$] the function $y\mapsto \varphi(x,y)$ is lower semicontinuous on 
$\RR^m$ for $\mathscr{H}^{n-1}$-a.e.\ $x\in M$.
\end{itemize}
A simple example is given by the function
\begin{eqnarray}\label{fi_ex}
\varphi(x,y)&:=& 
\begin{cases}
a + b|y| & \mbox{if } y\in \RR^m\smallsetminus \{0\},\\
0 & \mbox{if } y=0,
\end{cases}
\end{eqnarray}
where $a\geq 0$ and $b\geq 0$ are real constants.
The constant $a$ plays the role of an activation energy; if $b>0$, there is also an energy term proportional to the
amplitude of the crack opening. The classical Griffith's model corresponds to the case $a>0$ and $b=0$.

Let $L^0(M)$ be the set of extended real valued measurable functions on $M$ and let $L^0(M)^+$ be the set of functions $w\in L^0(M)$ 
such that $w\geq 0$ $\mathscr{H}^{n-1}$-a.e.\ on $M$.

We introduce the function $\phi\colon L^p(M;\RR^m) \to L^0(M)^+$ defined by
\begin{equation*}
\phi(w)(x):= \varphi(x, w(x)), 
\end{equation*}
for every $w \in L^p(M;\RR^m)$ and for $\mathscr{H}^{n-1}$-a.e.\  $x\in M$. 

Given an arbitrary family $(w_i)_{i\in I}$ in $L^0(M)^+$ the essential supremum 
\begin{equation*}
w = \esssup_{i\in I} \,w_i
\end{equation*}
of the family is defined as the unique  
(up to $\mathscr{H}^{n-1}$-equivalence) function in $L^0(M)^+$ such that
\begin{itemize}
\item[$\bullet$] $w \geq w_i$ $\mathscr{H}^{n-1}$-a.e.\  on $M$ for all $i\in I$;
\item[$\bullet$] if $z\in L^0(M)^+$ and $z\geq w_i$ $\mathscr{H}^{n-1}$-a.e.\ on $M$, then
$z\geq w$ $\mathscr{H}^{n-1}$-a.e.\ on $M$.
\end{itemize}
For the existence of such a function see, for instance, \cite[Proposition VI-{\bf 1}-1]{N-75}.

Suppose now that the deformation $u$ depends on time, i.e., we have a map $t\mapsto u(t)$ 
from $[0,T]$ into $W^{1,p}(\Omega\smallsetminus M;\RR^m)$.
If no crack is present until time $0$ and 
\begin{equation*}
\phi([u(s)])\leq \phi([u(t)]) \qquad \mathscr{H}^{n-1}\hbox{-a.e.\ on }M
\end{equation*}
for every $s\in[0, t]$,
then the energy dissipated in the crack process in the time interval $[0,t]$ is given, in our model, by
\begin{equation*}
\int_M \phi([u(t)])\,d\mathscr{H}^{n-1}.
\end{equation*}
This happens for instance when $s\mapsto\phi([u(s)])$ is 
monotonically increasing $\mathscr{H}^{n-1}$-a.e.\ on $M$.

In the general case, the irreversibility of the fracture process leads to
introduce an auxiliary function $t\mapsto\beta(t)$ from $[0,T]$ to $L^1(M)$, 
which takes into account the history of the system up to time $t$. 
We assume that for every $0\leq t_1\leq t_2\leq T$ we have 
\begin{equation}\label{b:cond}
\beta(t_2)= \beta(t_1)\vee\esssup_{t_1\leq s \leq t_2} \phi([u(s)])\quad \text{$\mathscr{H}^{n-1}$-a.e.\ on $M$,}
\end{equation}
so that
\begin{equation*}
\beta(t_2) - \beta(t_1)= \esssup_{t_1\leq s\leq t_2} (\phi([u(s)])-\beta(t_1))^+ \quad \text{$\mathscr{H}^{n-1}$-a.e.\ on $M$},
\end{equation*}
where for every $a\in \RR$, $a^+:= a\vee 0$ denotes the positive part of $a$.

In particular 
\begin{itemize}
\item[$\bullet$] $t\mapsto \beta(t)$ is increasing, i.e., $\beta(t_1)\leq \beta(t_2)$ $\mathscr{H}^{n-1}$-a.e.\ 
on $M$ for $0\leq t_1\leq t_2\leq T$; 
\item[$\bullet$] $\beta(t) \geq \phi([u(t)])$ $\mathscr{H}^{n-1}$-a.e.\ on $M$ for every $t\in [0,T]$.
\end{itemize}

In our model the energy dissipated in the time interval $[t_1,t_2]$ is given by
\begin{equation*}
\|\beta(t_2)-\beta(t_1)\|_{1,M}:=\int_M(\beta(t_2)-\beta(t_1))\,d\mathscr{H}^{n-1}  .
\end{equation*}
According to this assumption there is no dissipation in
the intervals $[t_1,t_2]$ where $\phi([u(s)])\leq \beta(t_1)$
$\mathscr{H}^{n-1}$-a.e.\ on $M$ for every $s\in [t_1,t_2]$,
while the dissipation is given by
\begin{equation*}
\int_M\big(\phi([u(t_2)])-\phi([u(t_1)])\big)\,d\mathscr{H}^{n-1}
\end{equation*}
whenever $\beta(t_1)\leq \phi([u(s)])\leq \phi([u(t_2)])$ for every $s\in [t_1,t_2]$.

It follows from (\ref{b:cond}) that $\beta(t)$ is uniquely determined by $\beta(0)$ and 
by the history of the deformation $s\mapsto u(s)$ in the interval $[0,t]$.
Since it is difficult to deal with (\ref{b:cond}) directly, 
we prefer to define the notion of quasistatic evolution by considering a more general internal variable 
$t\mapsto \gamma(t)$ which is assumed to satisfy the following weaker conditions:
\begin{itemize}
\item[$\bullet$] $t\mapsto \gamma(t)$ is increasing, i.e., $\gamma(t_1)\leq \gamma(t_2)$ $\mathscr{H}^{n-1}$-a.e.\ 
on $M$ for $0\leq t_1\leq t_2\leq T$; 
\item[$\bullet$] $\gamma(t) \geq \phi([u(t)])$ $\mathscr{H}^{n-1}$-a.e.\ on $M$ for every $t\in [0,T]$.
\end{itemize}
We do not assume from the beginning that $t\mapsto \gamma(t)$ satisfies (\ref{b:cond}). 
This property will be a nontrivial consequence of the other conditions considered in the definition of 
quasistatic evolution (see Theorem \ref{tm:essup}).

Given functions $\psi \in W^{1,p}(\Omega;\RR^m)$ and $\gamma\in L^0(M)^+$, 
it is convenient to introduce the set $AD(\psi,\gamma)$ of {\em admissible deformations} 
with boundary value $\psi$ on $\partial_0\Omega$ and internal variable $\gamma$. It is defined by
\begin{equation*}
AD(\psi,\gamma):= \{ u\in W^{1,p}(\Omega \smallsetminus M;\RR^m): 
\phi([u])\leq \gamma \text{ on $M$, and } 
u = \psi \text{ on } \partial_0\Omega\},
\end{equation*}
where equalities and inequalities are considered $\mathscr{H}^{n-1}$-a.e., and the last equality 
refers to the traces of $u$ and $\psi$ on $\partial_0\Omega$.

An \emph{admissible configuration} with boundary value $\psi$ on $\partial_0\Omega$ is a pair $(u,\gamma)$, 
with $\gamma \in L^1(M)^+:=L^1(M)\cap L^0(M)^+$ and $u\in AD(\psi,\gamma)$.
\end{section}

\begin{section}{Definition and properties of quasistatic evolutions}\label{def_qe}
For every $t\in[0,T]$, the {\em total energy}
of an admissible configuration $(u,\gamma)$ at time $t$ is defined as
\begin{equation*}
\mathscr{E}(t)(u,\gamma):= \mathcal{W}(\nabla u) 
-\langle\mathscr{L}(t),u\rangle 
+\|\gamma\|_{1,M},
\end{equation*}
where $\|\cdot\|_{1,M}$ denotes the $L^1$-norm on $M$.

We now introduce the following definition in the spirit of Griffith's original theory on the crack propagation.
\begin{definition}\label{gs} 
A pair $(u,\gamma)\in W^{1,p}(\Omega \smallsetminus M;\RR^m)\times L^1(M)^+$ 
is {\em globally stable} at time $t\in [0,T]$ if $u\in AD(\psi(t),\gamma)$
and 
\begin{equation}\label{egs}
\mathscr{E}(t)(u,\gamma)\leq\mathscr{E}(t)(v,\delta)
\end{equation}
for every $\delta\geq \gamma$ and for every $v\in AD(\psi(t),\delta)$.
\end{definition}

In other words, the total energy of $(u,\gamma)$ at time $t$ cannot be reduced by increasing the
internal variable $\gamma$ or by choosing a new admissible deformation with the same boundary condition.

\begin{remark}\label{r:apriori}
For every $t\in[0,T]$ let $(u(t),\gamma(t))\in W^{1,p}(\Omega \smallsetminus M;\RR^m)\times L^1(M)^+$ be globally stable at time $t$.
By Definition \ref{gs} we can deduce an a priori estimate on $u(t)$. Indeed, by comparing $\mathscr{E}(t)(u(t),\gamma(t))$ 
with $\mathscr{E}(t)(\psi(t),\gamma(t))$, which is
bounded uniformly with respect to $t$, we get 
that $\mathcal{W}(\nabla u(t))-\langle\mathscr{L}(t),u(t)\rangle$ is bounded uniformly in time.
Next, by the assumption (\ref{whp}) on $\mathcal{W}$
and the boundedness of $\mathscr{L}(t)$ in $(W^{1,p}(\Omega \smallsetminus M;\RR^m))'$, we obtain that the $W^{1,p}$-norm 
of $u(t)$, $\|u(t)\|_{1,p}$, 
is bounded uniformly with respect to $t$.
Furthermore from this fact and by Definition~\ref{gs} we get that the crack term $\|\gamma(t)\|_{1,M}$ 
is bounded uniformly in time, too. 
\end{remark}

\begin{remark}\label{r:equiv}
Condition $(\ref{egs})$ is equivalent to 
\begin{equation*}
\mathscr{E}(t)(u,\gamma)\leq\mathscr{E}(t)(v,\gamma \vee \phi([v])),
\end{equation*}         
for every $v \in W^{1,p}(\Omega \smallsetminus M;\RR^m)$ such that $v=\psi(t)$
$\mathscr{H}^{n-1}$-a.e.\ on $\partial_0\Omega$. This is equivalent to 
\begin{equation}\label{eq:egs2}
\mathcal{W}(\nabla u) -\langle \mathscr{L}(t), u\rangle \leq
\mathcal{W}(\nabla v) -\langle \mathscr{L}(t), v\rangle +
\|(\phi([v])-\gamma)^+\|_{1,M} 
\end{equation}     
for every $v \in W^{1,p}(\Omega \smallsetminus M;\RR^m)$ such that $v=\psi(t)$ $\mathscr{H}^{n-1}$-a.e.\ on $\partial_0\Omega$. 
This implies that if $(u,\gamma)$ is globally stable at time $t$ and $\tilde{\gamma}\in L^1(M)^+$ satisfies 
$\phi([u])\leq \tilde{\gamma}\leq \gamma$ $\mathscr{H}^{n-1}$-a.e.\ on $M$, 
then $(u,\tilde{\gamma})$ is globally stable at time~$t$.                                                             
\end{remark}

\begin{definition}\label{iqe}
An {\em irreversible quasistatic evolution} of minimum energy configurations is a 
function $t \mapsto (u(t),\gamma(t))$ from $[0,T]$ into $W^{1,p}(\Omega\smallsetminus M; \RR^m)\times L^1(M)^+$ 
which satisfies the following conditions:
\begin{itemize}
\item[(a)] {\em global stability}: for every $t\in[0,T]$
the pair $(u(t),\gamma(t))$ is globally stable at time $t$;
\item[(b)] {\em irreversibility}: $\gamma(s) \leq \gamma(t)$ $\mathscr{H}^{n-1}$-a.e.\ on $M$ for every $0\leq s\leq t\leq T$;
\item[(c)] {\em energy balance}: the function $t\mapsto \mathscr{E}(t)(u(t),\gamma(t))$ is absolutely continuous on $[0,T]$ and 
\begin{equation*}
\qquad\frac{d}{dt}(\mathscr{E}(t)(u(t),\gamma(t)))= \langle \partial\mathcal{W}(\nabla u(t)),\nabla \dot{\psi}(t)\rangle
-\langle\mathscr{L}(t),\dot{\psi}(t)\rangle-\langle\dot{\mathscr{L}}(t),u(t)\rangle\, ,
\end{equation*}
for a.e.\ $t\in[0,T]$. 
\end{itemize}
\end{definition}
\begin{remark}\label{r:c-c'}
Condition (c) is equivalent to the following one:
\begin{itemize}
\item[(c')] {\em energy balance in integral form}: the function 
$t\mapsto\langle \partial\mathcal{W}(\nabla u(t)),\nabla \dot{\psi}(t)\rangle
-\langle\dot{\mathscr{L}}(t),u(t)\rangle $  belongs to $L^1([0,T])$ and 
\begin{eqnarray*}
&\mathscr{E}(t)(u(t),\gamma(t)) - \mathscr{E}(0)(u(0),\gamma(0)) =\\
&= \displaystyle \int_0^t 
\Big(\langle \partial\mathcal{W}(\nabla u(s)),\nabla \dot{\psi}(s)\rangle
-\langle\mathscr{L}(s),\dot{\psi}(s)\rangle-\langle\dot{\mathscr{L}}(s),u(s)\rangle\Big)\, ds
\end{eqnarray*}
for every $t\in [0,T]$.
\end{itemize}
This can be written in the form 
\begin{eqnarray}\label{cc'}
&\mathcal{W}(\nabla u(t)) - \mathcal{W}(\nabla u(0)) +
\|\gamma(t) - \gamma(0)\|_{1,M}=\nonumber\\
& \displaystyle =\int_0^t\Big(\langle \partial\mathcal{W}(\nabla u(s)),\nabla \dot{\psi}(s)\rangle
-\langle\mathscr{L}(s),\dot{\psi}(s)\rangle\Big) \, ds +\\
&\displaystyle{} +\langle{\mathscr{L}}(t),u(t)\rangle - \langle{\mathscr{L}}(0),u(0)\rangle
- \int_0^t \langle\dot{\mathscr{L}}(s),u(s)\rangle\, ds,\nonumber
\end{eqnarray}
for every $t\in [0,T]$. 
The first line is the increment in stored energy plus a term which will be interpreted as the energy dissipated 
by the crack process in the time interval $[0,t]$, as we shall see in Remark~\ref{bg}. 
Using the divergence theorem we can show that the second line represents the work done in the same time interval by the forces 
which act on $\partial_0 \Omega$ to produce the imposed deformation. 
The third line represents the work done by the imposed forces in the interval $[0,t]$; 
this follows from an integration by parts when $t\mapsto u(t)$ is regular enough, and can be obtained 
by approximation in the other cases.

If $t\mapsto (u(t),\gamma(t))$ satisfies condition (a), then $(u(t),\gamma(t))$ is bounded in 
$W^{1,p}(\Omega\smallsetminus M;\RR^m)\times L^1(M)^+$ by Remark~\ref{r:apriori}. 
Therefore in condition (c') it is enough to assume that 
$t\mapsto\langle \partial\mathcal{W}(\nabla u(t)),\nabla \dot{\psi}(t)\rangle
-\langle\dot{\mathscr{L}}(t),u(t)\rangle $ is measurable.
\end{remark}

In the following theorem we prove one inequality of the energy balance.
\begin{theorem}\label{th:en>}
Let $t\mapsto (u(t),\gamma(t))$ be a function from $[0,T]$ into $W^{1,p}(\Omega\smallsetminus M;\RR^m)\times L^1(M)^+$ 
which satisfies the global stability condition $(a)$ and the irreversibility condition $(b)$ 
of Definition~\ref{iqe}. Assume that  $t\mapsto\langle \partial\mathcal{W}(\nabla u(t)),\nabla \dot{\psi}(t)\rangle
-\langle\dot{\mathscr{L}}(t),u(t)\rangle $ is measurable. Then
\begin{eqnarray*}
&\mathscr{E}(t)(u(t),\gamma(t)) - \mathscr{E}(0)(u(0),\gamma(0)) \geq\nonumber\\
&\geq \displaystyle \int_0^t 
\Big(\langle \partial\mathcal{W}(\nabla u(s)),\nabla \dot{\psi}(s)\rangle
-\langle\mathscr{L}(s),\dot{\psi}(s)\rangle-\langle\dot{\mathscr{L}}(s),u(s)\rangle\Big)\, ds
\end{eqnarray*}
for every $t\in [0,T]$.
\end{theorem}
\begin{proof}
We note that  $t\mapsto\langle \partial\mathcal{W}(\nabla u(t)),\nabla \dot{\psi}(t)\rangle
-\langle\dot{\mathscr{L}}(t),u(t)\rangle $ belongs to $L^1([0,T])$ by the estimates of Remark~\ref{r:apriori}. 
The result can now be obtained arguing as in \cite{DM-F-T-03} (see the proof of Lemma~7.1 and the final part of the proof of 
Theorem~3.15).
\end{proof}
Now we prove that for a quasistatic evolution $t\mapsto (u(t),\gamma(t))$, the internal variable $t\mapsto \gamma(t)$ satisfies
a condition analogous to (\ref{b:cond}).
\begin{theorem}\label{tm:essup}
Let $t\mapsto (u(t),\gamma(t))$ be a quasistatic evolution. Then
\begin{equation}\label{esssup}
\gamma(t_2) = \gamma(t_1) \vee \esssup_{t_1\leq s\leq t_2} \phi([u(s)])
\quad \mathscr{H}^{n-1}\mbox{-a.e.\ on }M,
\end{equation}
 for every $0\leq t_1\leq t_2\leq T$.
\end{theorem}
\begin{proof}
It is enough to prove that
\begin{equation}\label{esssup0}
\gamma(t) = \gamma(0) \vee \esssup_{0\leq s\leq t} \phi([u(s)])
\quad \mathscr{H}^{n-1}\mbox{-a.e.\ on }M,
\end{equation}
for every $t\in[0,T]$.
Let $\tilde{\gamma}(t)$ be the right-hand side of (\ref{esssup0}).  
Since $t\mapsto \gamma(t)$ is increasing and $ \phi([u(t)])\leq \gamma(t)$ $\mathscr{H}^{n-1}$-a.e.\ on $M$ 
for every $t\in [0,T]$, it follows that
$\tilde{\gamma}(t)\leq \gamma(t)$ $\mathscr{H}^{n-1}$-a.e.\ on $M$ for every $t\in [0,T]$. 
As $ \phi([u(t)])\leq\tilde{\gamma}(t)$ $\mathscr{H}^{n-1}$-a.e.\ on $M$, by Remark~\ref{r:equiv} 
the pair $(u(t),\tilde{\gamma}(t))$ is globally stable at time $t$ for every $t\in [0,T]$. 
Since $t\mapsto \tilde{\gamma}(t)$ is increasing, we can apply Theorem~\ref{th:en>} and we obtain  
\begin{eqnarray*}
&\mathscr{E}(t)(u(t),\tilde{\gamma}(t)) - \mathscr{E}(0)(u(0),\gamma(0)) \geq\nonumber\\
&\geq \displaystyle \int_0^t 
\Big(\langle \partial\mathcal{W}(\nabla u(s)),\nabla \dot{\psi}(s)\rangle
-\langle\mathscr{L}(s),\dot{\psi}(s)\rangle-\langle\dot{\mathscr{L}}(s),u(s)\rangle\Big)\, ds
\end{eqnarray*}
for every $t\in [0,T]$. By the energy balance (c) it follows that 
$\mathscr{E}(t)(u(t),\tilde{\gamma}(t)) \geq \mathscr{E}(t)(u(t),\gamma(t))$, i.e.,
\begin{eqnarray*}
\mathcal{W}(\nabla u(t)) - \langle \mathscr{L}(t),u(t)\rangle + \|\tilde{\gamma}(t)\|_{1,M}
\geq \mathcal{W}(\nabla u(t)) - \langle \mathscr{L}(t),u(t)\rangle + \|\gamma(t)\|_{1,M},
\end{eqnarray*}
which implies $\|\tilde{\gamma}(t)\|_{1,M} \geq \|\gamma(t)\|_{1,M}$. 
As $\tilde{\gamma}(t)\leq \gamma(t)$ $\mathscr{H}^{n-1}$-a.e.\ on $M$, we deduce that 
$\tilde{\gamma}(t) = \gamma(t)$ $\mathscr{H}^{n-1}$-a.e.\ on $M$ for every $t\in [0,T]$, which concludes the proof.
\end{proof}
Theorem~\ref{tm:essup} can be used to explain the mechanical meaning of the internal variable $\gamma$ in the model 
case $\varphi(x,y):= |y|$. Indeed, if $t\mapsto (u(t),\gamma(t))$ is a quasistatic evolution with $\gamma(0)=0$ and 
$\varphi(x,y):= |y|$, then (\ref{esssup}) shows that $\gamma(t)(x)$ coincides with the maximum modulus of the 
amplitude of the opening reached by the crack at $x$ up to time $t$. 
\begin{remark}\label{bg}
As $t\mapsto \gamma(t)$ satisfies (\ref{b:cond}) by Theorem~\ref{tm:essup},
the mechanical interpretation given in Section \ref{setting} shows that
the term ${\|\gamma(t) - \gamma(0)\|_{1,M}}$ in (\ref{cc'}) represents the energy dissipated in the crack 
process in the time interval $[0,t]$.
\end{remark}
\begin{remark}\label{r:mielke}
In our model, the dissipation term in the energy functional comes from the expression $\|\gamma\vee\phi([v])-\gamma\|_{1,M}$ 
and is nonlinear in $\gamma$. 
This turns out to be the main mathematical difference between our model and the model considered by Mielke and Mainik and Mielke 
in \cite[Section~4.2]{Mie-Beijing} and \cite[Section~6.2]{M-M-03}, where the dissipation term is linear.
\end{remark}
We are now in a position to state our main result.
\begin{theorem}\label{mainthm}
Let $(u_0,\gamma_0) \in W^{1,p}(\Omega\smallsetminus M;\RR^m)\times L^1(M)^+$ be globally stable at time $t=0$. Then there exists
an irreversible quasistatic evolution $t\mapsto (u(t),\gamma(t))$ such that 
$(u(0),\gamma(0))=(u_0,\gamma_0)$.
\end{theorem}
\end{section}

\begin{section}{Some tools}\label{tools}

We introduce a notion of convergence for the functions $\gamma$, which is the counterpart of the notion of 
convergence of sets introduced in \cite{DM-F-T-03}.
The main property of this convergence is that, if $u_k$ 
converges weakly in $W^{1,p}(\Omega\smallsetminus M;\RR^m)$ to some function $u$ and 
$\phi([u_k])\leq \gamma_k$ ${\mathscr H}^{n-1}$-a.e.\ on $M$, then $\phi([u])\leq\gamma$ 
${\mathscr H}^{n-1}$-a.e.\ on $M$.

\begin{definition}\label{sigmap}
Let $\gamma_k, \gamma\in L^0(M)^+$. 
We say that $\gamma_k$ $\sigma^p_\varphi$-converges to $\gamma$ if the following two conditions are satisfied:
\begin{itemize}
\item[(a)] if $u_j \rightharpoonup u$ weakly in $W^{1,p}(\Omega\smallsetminus M;\RR^m)$ and 
$\phi([u_j])\leq \gamma_{k_j}$ ${\mathscr H}^{n-1}$-a.e.\ on $M$ for some sequence $k_j\to\infty$, then 
$\phi([u])\leq \gamma$ ${\mathscr H}^{n-1}$-a.e.\ on $M$;
\item[(b)] there exist a sequence  $u^i\in W^{1,p}(\Omega\smallsetminus M;\RR^m)$, with 
$\sup_{i} \phi([u^i])=\gamma$  ${\mathscr H}^{n-1}$-a.e.\ on $M$, and, for every $i$, a sequence 
$u^i_k\in W^{1,p}(\Omega\smallsetminus M;\RR^m)$, converging to $u^i$ weakly in 
$W^{1,p}(\Omega\smallsetminus M;\RR^m)$ as $k\to \infty$, such that $\phi([u^i_k])\leq\gamma_k$ 
${\mathscr H}^{n-1}$-a.e.\ on $M$ for every~$i$ and~$k$. 
\end{itemize}
\end{definition} 
Notice that we do not require any upper bound in $L^1(M)^+$ for the functions~$\gamma_k$.
\begin{remark}\label{rem_sp}
If $\gamma_k$ $\sigma^p_\varphi$-converges to $\gamma$, then in particular
there are functions $u^i_k$ and $u^i$ in $W^{1,p}(\Omega\smallsetminus M;\RR^m)$ such that 
condition (b) in Definition~$\ref{sigmap}$ holds.
We define for every $k$
\begin{equation*}
\gamma^i:=\sup_{j=1,\dots,i} \phi([u^j]) \quad \mbox{and}\quad
\gamma^i_k:=\sup_{j=1,\dots,i} \phi([u^j_k]). 
\end{equation*}
With this notation it turns out that
\begin{equation*}
\gamma = \lim_{i\to \infty} \gamma^i \quad \text{ and }
\quad \gamma_k \geq \sup_{i\in \NN} \gamma^i_k,
\end{equation*}
for every $k$.
\end{remark}
\begin{remark}\label{sp-point}
If $\gamma_k$ $\sigma^p_\varphi$-converges to $\gamma$, then 
\begin{equation*}
\gamma\leq \limsup_{k\to\infty} \gamma_k, \quad \mathscr{H}^{n-1}\mbox{-a.e.\  on } M,
\end{equation*}
as we can see by modifing the proof of Lemma~\ref{lm_lscsigma} below. 
Notice that the inequality can be strict, even when $\gamma_k$ converges pointwise to a function $\tilde{\gamma}$.
As an example, consider $n=2$, $m=1$, $p=2$, $\Omega={]{-2},2[}^2$ and $M=[0,1]\times \{0\}$.
Let $\gamma_k\in L^0(M)^+$ be defined as follows:
\begin{equation*}
\gamma_k(x):=
\begin{cases}
1 & \hbox{for } x\in\left[\frac{i}{k},\frac{i+1}{k}-\frac{1}{k^2}\right[;\\ \\
0 & \hbox{for } x\in\left[\frac{i+1}{k}-\frac{1}{k^2},\frac{i+1}{k}\right[;
\end{cases}
\quad \hbox{for } i=0,\dots,k-1.
\end{equation*}
It follows from homogenization theory (see \cite{Da-85}, \cite{Mu-85}, \cite{Pi-87})
that condition (a) in Definition~\ref{sigmap} is satisfied with 
$\gamma=0$, hence $\gamma_k$ $\sigma^2_\varphi$-converges to $0$. 
Furthermore $\gamma_k$ converge in measure to $1$, so up to a subsequence we have pointwise convergence
to $1=:\tilde{\gamma}>\gamma$.
\end{remark} 
We prove in the following lemma that the $L^1$-norm is lower semicontinuous with respect to $\sigma^p_\varphi$-convergence.
\begin{lemma}\label{lm_lscsigma}
Let $\gamma_k$, $\gamma\in L^0(M)^+$.
If $\gamma_k$ $\sigma^p_\varphi$-converges to $\gamma$
then
\begin{equation}\label{lsc_sp}
\|\gamma\|_{1,M}\leq 
\liminf_{k\to\infty} \|\gamma_k\|_{1,M}.
\end{equation} 
\end{lemma}
\begin{proof}
{}From the hypothesis it follows in particular that there are 
functions $u^i_k$ and $u^i$ in $W^{1,p}(\Omega\smallsetminus M;\RR^m)$ which satisfy condition (b) in Definition~\ref{sigmap}.
With notation from Remark~\ref{rem_sp}, let us prove that for every $i$
\begin{equation}\label{lsc_sp1}
\|\gamma^i\|_{1,M}\leq 
\liminf_{k\to\infty} \|\gamma^i_k\|_{1,M}.
\end{equation} 
Extracting a subsequence we may assume that $\liminf_{k}\|\gamma^i_k\|_{1,M}$ is a limit. 
As $[u^j_k]\to[u^j]$ strongly in $L^p(M;\RR^m)$ for $j=1,\dots,i$, we can extract a further subsequence such that 
$[u^j_k]\to[u^j]$ pointwise $\mathscr{H}^{n-1}$-a.e.\ on $M$ for $j=1,\dots,i$.
By the lower semicontinuity assumption $(\varphi_3)$ this implies 
\begin{equation*}
\gamma^i\leq  \liminf_{k\to\infty}\gamma^i_k \quad \mathscr{H}^{n-1}\mbox{-a.e.\ on } M.
\end{equation*}
By the Fatou lemma we obtain (\ref{lsc_sp1}), which yields
\begin{equation*}
\|\gamma^i\|_{1,M}\leq 
\liminf_{k\to\infty} \|\gamma_k\|_{1,M}.
\end{equation*}
We then pass to the limit as $i$ tends to infinity and obtain (\ref{lsc_sp}).
\end{proof}
We now prove a compactness result for the notion of $\sigma^p_\varphi$-convergence.
\begin{lemma}\label{comp_sp}
Every sequence in $L^0(M)^+$ has a $\sigma^p_\varphi$-convergent subsequence.
\end{lemma}
\begin{proof}
Let us denote the $L^p$-norm by $\|\cdot\|_p$.
Let $\gamma_k\in L^0(M)^+$, let $w_h\in L^{\infty}(\Omega\smallsetminus M;\RR^m)$ be dense in $L^p(\Omega\smallsetminus M;\RR^m)$, and, 
for every positive integers $l$, $h$, and $k$, let us consider the problem
\begin{equation}\label{yoshida}
\min \big\{\left\|\nabla u\right\|_p^p + \ell \left\|u-w_h\right\|_p^p \big\},
\end{equation}
where the minimum is taken over all functions $u\in W^{1,p}(\Omega\smallsetminus M;\RR^m)$ such that $\phi([u])\leq \gamma_k \;\;{\mathscr H}^{n-1}$-a.e.\ on $M$.

To prove that the minimum is achieved, we take a minimizing sequence and we easily obtain that it is bounded in
$W^{1,p}(\Omega\smallsetminus M;\RR^m)$.
Then, up to a subsequence, we can pass to the limit and by using our lower semicontinuity assumption $(\varphi_3)$
we can prove that the limit function is actually a solution to the minimum problem (\ref{yoshida}).
This solution, which is unique by strict convexity, will be denoted by $u^{\ell,h}_k$.
Notice that this function is bounded in $W^{1,p}(\Omega\smallsetminus M;\RR^m)$ uniformly with respect to $k$,
thus, up to a subsequence, we can pass to the limit in $k$ 
and get that there is a function $u^{\ell,h}$ such that
$u^{\ell,h}_k\rightharpoonup u^{\ell,h}$ weakly in $W^{1,p}(\Omega\smallsetminus M;\RR^m)$.
Further we define
\begin{equation}\label{def:g}
\gamma:= \sup_{\ell,h\in \NN} \phi([u^{\ell,h}]) \quad {\mathscr H}^{n-1}\text{-a.e.\ on $M$.}
\end{equation}
In this way point (b) of Definition~\ref{sigmap} is automatically satisfied. 

We need to prove point (a).
To this aim, let
$v_j \rightharpoonup v$ weakly in $W^{1,p}(\Omega\smallsetminus M;\RR^m)$ be such that
$\phi([v_j])\leq \gamma_{k_j}$ ${\mathscr H}^{n-1}$-a.e.\ on $M$ for some sequence $k_j\to\infty$.
We want to prove that
$\phi([v])\leq \gamma$ ${\mathscr H}^{n-1}$-a.e.\ on $M$. By density there is a subsequence of $w_h$, say $w_{h_i}$,
which converges strongly to $v$ in $L^p(\Omega\smallsetminus M;\RR^m)$. Let $\ell_i\to+\infty$ be such that
$\ell_i\|v - w_{h_i}\|_p^p\to 0$ as $i$ tends to infinity.
By the minimality of $u^{\ell_i,h_i}_{k_j}$, we have 
\begin{equation*}
\|\nabla u^{\ell_i,h_i}_{k_j}\|_p^p+\ell_i\|u^{\ell_i,h_i}_{k_j}-w_{h_i}\|_p^p
\leq \|\nabla v_j \|_p^p + \ell_i\|v_j - w_{h_i}\|_p^p.
\end{equation*}
Then $u^{\ell_i,h_i}_{k_j}$ is bounded in $W^{1,p}(\Omega\smallsetminus M;\RR^m)$ uniformly with respect to $j$, 
and passing to the limit as $j$ tends to infinity we get 
\begin{equation*}
\|\nabla u^{\ell_i,h_i}\|_p^p+\ell_i\|u^{\ell_i,h_i}-w_{h_i}\|_p^p
\leq \sup_{j\in\NN}\|\nabla v_{j}\|_p^p+\ell_i\|v-w_{h_i}\|_p^p.
\end{equation*}
Since $\ell_i\|v - w_{h_i}\|_p^p\to 0$ as $i$ tends to infinity, 
this inequality ensures that $\nabla u^{\ell_i,h_i}$ is bounded
in $L^p(\Omega\smallsetminus M;\RR^m)$ uniformly with respect to $i$, and $u^{\ell_i,h_i}_{k_j}-w_{h_i} \to 0$ 
strongly in $L^p(\Omega\smallsetminus M;\RR^m)$. As $w_{h_i}\to v$ strongly in $L^p(\Omega\smallsetminus M;\RR^m)$, 
we deduce that $u^{\ell_i,h_i}$ converges weakly to $v$
in $W^{1,p}(\Omega\smallsetminus M;\RR^m)$. Then $[u^{\ell_i,h_i}]$ converges strongly to $[v]$ in $L^p(M;\RR^m)$.
Passing to a subsequence, we may also obtain pointwise convergence $\mathscr{H}^{n-1}$-a.e.\ on $M$.
By (\ref{def:g}) we have $\phi([u^{\ell_i,h_i}])\leq \gamma$ $\mathscr{H}^{n-1}$-a.e.\ on $M$,
so that the lower semicontinuity assumption $(\varphi_3)$ yields $\phi([v])\leq \gamma$ $\mathscr{H}^{n-1}$-a.e.\ on $M$,
which is precisely the conclusion to point (a) in the definition of $\sigma^p_\varphi$-convergence.
\end{proof}
We shall use the following Helly-type compactness result. We recall that a function 
$t\mapsto \gamma(t)$  from $[0,T]$ into $L^0(M)^+$ is said to be increasing if 
$\gamma(s)\leq \gamma(t)$ $\mathscr{H}^{n-1}$-a.e.\ on $M$, whenever $0\leq s\leq t\leq T$. 
\begin{lemma}\label{lm:helly}
Let $t\mapsto \gamma_k(t)$ be a sequence of increasing functions from $[0,T]$ into $L^0(M)^+$. 
Then there exist a subsequence $\gamma_{k_j}$, independent of $t$, and an increasing function 
$t\mapsto \gamma(t)$ from $[0,T]$ into $L^0(M)^+$, 
such that $\gamma_{k_j}(t)$  $\sigma^p_\varphi$-converges 
to $\gamma(t)$ for every $t\in [0,T]$.
\end{lemma}
\begin{proof}
Let $D$ be a countable dense subset of $[0,T]$ containing $0$ and $T$. 
By Lemma~\ref{comp_sp}, using a diagonal argument, we can extract a
 subsequence, still named $\gamma_{k}(t)$,  and an increasing function
$t\mapsto \gamma(t)$ from $D$ into $L^0(M)^+$, such that
$\gamma_{k}(t)$ $\sigma^p_\varphi$-converges to $\gamma(t)$ for every $t\in D$.

Let us define
\begin{equation*}
\gamma(t+) := \inf_{s\geq t, \, s\in D} \gamma(s)\qquad \hbox{and} \qquad
\gamma(t-) := \sup_{s\leq t, \, s\in D} \gamma(s),
\end{equation*}
for every $t\in[0,T]$. It is easy to prove that:
\begin{itemize}
\item[(1)] $\gamma(t-)= \gamma(t) = \gamma(t+)$ for every $t\in D$;
\item[(2)] $\gamma(t-)\leq \gamma(t+)$ for every $t\in [0,T]$;
\item[(3)] if $s < t$, then $\gamma(s+) \leq \gamma(t-)$.
\end{itemize}
Define $E:= \{t\in [0,T]:\gamma(t+)= \gamma(t-)\,\mathscr{H}^{n-1}\hbox{-a.e.\ in }M\}$ and 
$\gamma(t):=\gamma(t-)= \gamma(t+)$ for every $t\in E$. Note that by (1) $D$ is contained in $E$ 
and the definition of $\gamma(t)$ agrees with the original one on $D$. Then the definition of 
$\sigma^p_\varphi$-convergence and the monotonicity condition imply that $\gamma_k(t)$ $\sigma^p_\varphi$-converges to $\gamma(t)$
for every $t\in E$.

Let us show now that the set $E^c:= [0,T]\smallsetminus E$ is at most countable.
For every pair of positive integers $i,k$ we set 
$A_{i,k} := \{t\in [0,T]: \|(\gamma(t+)\land k)- (\gamma(t-)\land k)\|_{1,M}> 1/ i\}$, so that
we have $E^c$ is the union of the sets $A_{i,k}$. Therefore it is enough to show that each set $A_{i,k}$ is finite.
Let $t_1<\dots <t_r \in A_{i,k}$.
Since, by (3), $(\gamma(t_{j-1}+)\land k)\leq ( \gamma(t_{j}-)\land k)$ for $j=2, \dots, r$, we get
\begin{equation*}
\frac r i \leq \sum_{j=1}^r  \|(\gamma(t_j+)\land k)-( \gamma(t_j-)\land k)\|_{1,M} \leq  
\|\gamma(t_r+)\land k\|_{1,M}\leq k \mathscr{H}^{n-1}(M) ,
\end{equation*}
so that $r\leq i k \mathscr{H}^{n-1}(M)$, which implies that 
$A_{i,k}$ is finite. It follows that $E^c$ is at most countable,
thus we can conclude the proof of the lemma by applying again the compactness 
Lemma~\ref{comp_sp}  for every $t\in E^c$, together with a diagonal argument.
\end{proof}
The following result plays a crucial role in the proof of point (a) in the Definition~\ref{iqe} of quasistatic evolution.
\begin{lemma}\label{lm_limsup}
Let $\gamma_k$, $\gamma\in L^0(M)^+$. 
Assume that $\gamma_k$ $\sigma^p_\varphi$-converges to $\gamma$.
Then for any $v\in W^{1,p}(\Omega\smallsetminus M;\RR^m)$ with $\phi([v])\in L^1(M)^+$
the following inequality holds true:
\begin{equation}\label{eq_limsup}
\limsup_{k\to\infty}  \|(\phi([v])-\gamma_k)^+\|_{1,M}\leq\|(\phi([v])-\gamma)^+\|_{1,M}.
\end{equation}
\end{lemma}
\begin{proof}
It is not restrictive to assume that the $\limsup$ is a limit.
Let $u^i$ and $u^i_k$ be the functions considered in point (b) of Definition~\ref{sigmap}. During the proof
we shall use the notation introduced in Remark~\ref{rem_sp}. As $\gamma_k^i\leq \gamma_k$ 
$\mathscr{H}^{n-1}$-a.e.\ on $M$, we have
\begin{equation*}
(\phi([v])-\gamma_k)^+ \leq (\phi([v])-\gamma_k^i)^+,
\end{equation*}
hence
\begin{equation}\label{keyeq0}
\lim_{k\to\infty}  \|(\phi([v])-\gamma_k)^+\|_{1,M}   \leq
\liminf_{k\to\infty}  \|(\phi([v])-\gamma_k^i)^+\|_{1,M}\, .
\end{equation}
Passing to a subsequence, we may assume that 
$[u^i_k]$ converges to 
$[u^i]$ $\mathscr{H}^{n-1}$-a.e.\ on $M$. 
By the lower semicontinuity assumption $(\varphi_3)$ we obtain
\begin{equation*}
\gamma^i\leq \liminf_{k\to\infty}\gamma^i_k \quad\mbox{$\mathscr{H}^{n-1}$-a.e.\ on $M$},
\end{equation*}
so that Fatou Lemma gives 
\begin{equation*}
\limsup_{k\to\infty} \|(\phi([v])-\gamma_k^i)^+\|_{1,M} \leq
\|(\phi([v])-\gamma^i)^+\|_{1,M},
\end{equation*}
which, together with (\ref{keyeq0}), yields 
\begin{equation*}
\lim_{k\to\infty}  \|(\phi([v])-\gamma_k)^+\|_{1,M} \leq
\|(\phi([v])-\gamma^i)^+\|_{1,M}.
\end{equation*}
As $\gamma^i\to\gamma$ $\mathscr{H}^{n-1}$-a.e.\ on $M$,  
inequality~(\ref{eq_limsup}) can be obtained by passing to the limit as $i\to\infty$.
\end{proof}
\begin{remark}\label{notw*}
The conclusion of Lemma~\ref{lm_limsup} does not hold, in general, when 
$\gamma_k, \gamma\in L^{\infty}(M)^+$ and $\gamma_k \rightharpoonup \gamma$ 
weakly* in $L^{\infty}(M)$. Consider, for instance, the case $n=2$, $m=1$, $\Omega={]{-4},4[}^2$, $M=[-\pi,\pi]\times \{0\}$,
and define $\gamma_k(x):= 1 + \sin(kx_1)$, where $x_1$ denotes the first coordinate of $x$. 
Then, $\gamma_k$ converges to $\gamma(x):=1$ weakly* in $L^{\infty}(M)$, but (\ref{eq_limsup}) is not satisfied for 
$\phi([v]) = 1$, since in this case  $ \|(\phi([v])-\gamma_k)^+\|_{1,M} =2$ for every $k$, while  $ \|(\phi([v])-\gamma)^+\|_{1,M}=0$. 
\end{remark}
\end{section}

\begin{section}{The discrete-time problems and proof of the main result}\label{mainresult}

In this section we prove Theorem~\ref{mainthm} by a discrete-time approximation.
We fix a sequence of subdivisions $(t_k^i)_{0\leq i \leq k}$ 
of the interval $[0,T]$, with
\begin{eqnarray}\label{tk1}
0 = t_k^0 < t_k^1 < \dots <t_k^{k-1} < t_k^k = T,\\
\label{tk2}
\lim_{k \to \infty} \max_{1 \leq i \leq k} (t_k^i-t_k^{i-1})=0. 
\end{eqnarray}
For $i=1, \dots, k$ we set $\mathscr{L}^i_k = \mathscr{L}(t_k^i)$, $\psi^i_k=\psi(t^i_k) $, 
$\mathscr{E}^i_k=\mathscr{E}(t^i_k)$.

For every $k\in \NN$ we define $u_k^i$ and $\gamma_k^i$ by induction as follows. 
Let $(u_0,\gamma_0)$ be a minimum energy configuration at time $t=0$. We set
$(u_k^0, \gamma_k^0):= (u_0,\gamma_0)$ and define $(u_k^i, \gamma_k^i)$ as
a solution of the minimum problem
\begin{equation}\label{Pik}
\min \big\{ 
\mathscr{E}^i_k( u,\gamma) : \gamma \in L^1(M)^+,   \; \gamma \geq \gamma^{i-1}_{k}, \; u\in AD(\psi^i_k,\gamma)\big\}\,,
\end{equation}
where the inequality means that $\gamma \geq \gamma^{i-1}_{k}$ $\mathscr{H}^{n-1}$-a.e.\ on $M$.
\begin{remark}\label{Emindiscr}
Consider the minimum problem
\begin{equation}\label{Pik2}
\min \Big\{ 
\mathcal{W}(\nabla u)
- \langle\mathscr{L}^i_k,u\rangle + \|\phi([u]) \vee \gamma^{i-1}_k \|_{1,M}:
u = \psi_{k}^{i}
\text{ on } \partial_0\Omega
\Big\},
\end{equation} 
where $u$ is assumed to belong to  $W^{1,p}(\Omega \smallsetminus M;\RR^m)$.
Then the following two conditions are equivalent:
\begin{itemize}
\item[(a)] the pair $(u^i_k,\gamma^i_k)$ is a solution to (\ref{Pik});
\item[(b)] $u^i_k$ is a solution to (\ref{Pik2}) and 
$\gamma^i_k := \gamma^{i-1}_k \vee \phi([u^i_k])$  $\mathscr{H}^{n-1}$-a.e.\ on $M$.
\end{itemize}
The existence of a solution of (\ref{Pik}) (or equivalently (\ref{Pik2})) 
can be easily obtained by using the direct methods of the calculus of variations. 
The compactness of a minimizing sequence follows from (\ref{whp}) and positiveness of $\varphi$. 
The lower semicontinuity follows from $(W_1)$, $(W_2)$, $(\varphi_3)$, and from the compactness of the trace operator.
\end{remark}
For every $t \in [0,T]$ we define
\begin{equation}\label{kt}
\begin{array}{c}
\tau_k(t)=t_k^{i}, \; u_k(t)=u_k^{i}, \; \gamma_k(t)=\gamma_k^{i},\;\psi_k(t)=\psi(t_k^{i}),\\[2mm]
\mathscr{L}_k(t)=\mathscr{L}(t_k^{i}), \;\mathscr{E}_k(t)=\mathscr{E}(t_k^{i}),
\end{array}
\end{equation}
where $i$ is the greatest integer such that $t_k^{i} \leq t$.
Note that $u_k(t)=u_k(\tau_k(t))$, $\gamma_k(t)=\gamma_k(\tau_k(t))$, $\psi_k(t)=\psi(\tau_k(t))$, 
$\mathscr{L}_k(t)=\mathscr{L}(\tau_k(t))$ and  $\mathscr{E}_k(t)=\mathscr{E}(\tau_k(t))$.

\begin{remark}\label{bounds2}
Since $\psi^i_k\in AD(\psi^i_k, \gamma^{i-1}_k)$, then by Remark~\ref{r:apriori} we deduce that the $L^p$-norms $\|\nabla u^i_k\|_p$ 
and $\|u^i_k\|_p$ are bounded uniformly with respect to $i$ and $k$. Passing to the piecewise constant functions $t\mapsto \nabla u_k(t)$
and $t\mapsto u_k(t)$, we have that there exists a positive constant $C$ such that
\begin{equation}\label{bd}
\|\nabla u_k(t)\|_{p}\leq C \quad \mbox{and} \quad\|u_k(t)\|_{p}\leq C 
\end{equation}
for every $k$ and for every $t\in [0,T]$. Since $\mathscr{E}_k(t)(u_k(t),\gamma_k(t))$ is bounded 
uniformly with respect to $k$, we get also that 
\begin{equation}\label{bdg}
\|\gamma_k(t)\|_{1,M}\leq C,
\end{equation}
for every $k$ and for every $t\in [0,T]$.
\end{remark}
We introduce now a sequence of functions which play an important role in our
estimates. For a.e.\ $t\in[0,T]$ we set
\begin{equation}\label{teta_k}
\theta_k(t) := \langle \partial\mathcal{W}(\nabla u_k(t)),\nabla \dot{\psi}(t)\rangle  -
\langle\mathscr{L}_k(t),\dot{\psi}(t)\rangle- \langle\dot{\mathscr{L}}(t),u_k(t)\rangle.
\end{equation}

In the following lemma we present the main energy estimate for the discrete process.
\begin{lemma}\label{energyestim}
There exists a sequence $R_k \to 0$ such that
\begin{equation*}
\mathscr{E}(\tau_k(t))(u_k(t), \gamma_k(t))\leq \mathscr{E}(0)(u_0,\gamma_0)+
\int_0^{\tau_k(t)} \theta_k(s) \, ds + R_k,
\end{equation*}
for every $k$ and for every $t \in [0,T]$.
\end{lemma}
\begin{proof} 
We need to prove that there exists a sequence $R_k\to 0$ such that
\begin{equation*}
\mathscr{E}^i_k(u_k^i, \gamma_k^i)\leq \mathscr{E}(0)(u_0,\gamma_0)+
\int_0^{t_k^i} \theta_k(s) \, ds + R_k,
\end{equation*}
for any $k$ and for any $i=1,\dots,k$.

Let us fix $j$ and $k$ with $1\leq j\leq k$. Since $u^{j-1}_k=\psi^{j-1}_k$
on $\partial_0 \Omega$, and $[u^{j-1}_k+\psi^j_k-\psi^{j-1}_k]= [u^{j-1}_k]$ $\mathscr{H}^{n-1}$-a.e.\ on $M$,
the function $u^{j-1}_k+\psi^j_k-\psi^{j-1}_k$ belongs to $AD(\psi^j_k,\gamma^{j-1}_k)$, hence
$\mathscr{E}^j_k(u^j_k,\gamma^j_k)\leq \mathscr{E}^j_k
(u^{j-1}_k+\psi^j_k-\psi^{j-1}_k, \gamma_k^{j-1})$. 
The proof now can be concluded arguing as in the proof of \cite[Lemma 6.1]{DM-F-T-03}.
\end{proof}
We are now in a position to prove our main result.

\begin{proof}[Proof of Theorem~\ref{mainthm}]
Let $(t^i_k)$, $0\leq i\leq k$, be a sequence of subdivisions of the interval $[0,T]$ satisfying
(\ref{tk1}) and (\ref{tk2}). For any $k$ consider the pairs $(u^i_k,\gamma^i_k)$ inductively defined as solutions
of the discrete problems (\ref{Pik}) for $i=1,\dots,k$ with the initial condition $(u^0_k,\gamma^0_k)=(u_0,\gamma_0)$. Let 
$\tau_k(t)$, $u_k(t)$, $\gamma_k(t)$, and $\psi_k(t)$
be defined by (\ref{kt}) for any $t\in[0,T]$. 
By Lemma~\ref{lm:helly} there exists a subsequence of $\gamma_k(t)$, independent of $t$, which $\sigma^p_\varphi$-converges 
to $\gamma_\infty(t)\in L^0(M)^+$, for every $t\in [0,T]$. By (\ref{bdg}) and Lemma~\ref{lm_lscsigma} we have 
$\gamma_\infty(t)\in L^1(M)^+ $.

Let $\theta_k(t)$ be defined by (\ref{teta_k}) for a.e.\ $t$ and let
\begin{equation*}
\theta_{\infty}(t):= \limsup_{k\to\infty} \theta_k(t).
\end{equation*} 
By (\ref{dewhp}) and (\ref{bd}) we deduce that
\begin{equation*}
|\theta_k(t)|\leq \alpha_2 (C^{p-1}+1)\|\nabla \dot{\psi}(t)\|_{p}+
\|\mathscr{L}_k(t)\|_{*}\|\dot{\psi}(t)\|_{1,p}+
C\|\dot{\mathscr{L}}(t)\|_{*},
\end{equation*}
where $\|\cdot\|_*$ is the norm in the dual space of $W^{1,p}(\Omega\smallsetminus M;\RR^m)$. 
Since the right-hand side of previous formula belongs to $L^1([0,T])$,
we deduce that $\theta_{\infty}$ belongs to $L^1([0,T])$, too,
and using the Fatou lemma  we get
\begin{equation}\label{fatoulimsup_tk}
\limsup_{k\to \infty} \int_0^{\tau_k(t)}\theta_k(s)\,ds
\leq \int_0^t \theta_{\infty}(s)\,ds.
\end{equation}
For a.e.\ $t\in [0,T]$ we can extract a subsequence $\theta_{k_j}$ of $\theta_k$,
depending on $t$, such that
\begin{equation*}
\theta_{\infty}(t)=\lim_{j \to \infty} \theta_{k_j}(t).
\end{equation*}
By (\ref{bd}) the sequence $u_{k_j}(t)$ is bounded in $W^{1,p}(\Omega\smallsetminus M;\RR^m)$, 
therefore we can extract a further subsequence, still denoted by $u_{k_j}(t)$, which converges weakly in 
$W^{1,p}(\Omega\smallsetminus M;\RR^m)$ to a function $u_\infty(t)$.

Since $\phi([u_{k_j}(t)])\leq \gamma_{k_j}(t)$ $\mathscr{H}^{n-1}$-a.e.\ on $M$, by point (a) in Definition~\ref{sigmap} 
we have $\phi([u_{\infty}(t)]) \leq \gamma_{\infty}(t)$
$\mathscr{H}^{n-1}$-a.e.\ on $M$. 
On the other hand, as $u_{k_j}(t)=\psi_{k_j}(t)$ $\mathscr{H}^{n-1}$-a.e.\
on $\partial_0\Omega$, we have also $u_{\infty}(t)=\psi(t)$ $\mathscr{H}^{n-1}$-a.e.\
on $\partial_0\Omega$, so that $u_\infty(t) \in AD(\psi(t),\gamma_{\infty}(t))$ for every $t\in [0,T]$.

The next step is to prove that the pair $(u_{\infty}(t),\gamma_{\infty}(t))$
satisfies property (a) of Definition~\ref{iqe}. To this aim, let $\gamma\in L^1(M)^+$,
$\gamma\geq\gamma_{\infty}(t)$ and $v\in AD(\psi(t),\gamma)$.
By the minimality of the incremental solutions $(u_k(t),\gamma_k(t))$, we have that
$\mathscr{E}_k(t)\left(u_{k}(t),\gamma_{k}(t)\right) \leq
\mathscr{E}_k(t)(v_k,\gamma_{k}(t)\vee\phi([v]))$, where $v_k:= v+\psi_k(t)-\psi(t)$. 
Since the functional $u\mapsto \mathcal{W}(\nabla u)$ is weakly lower semicontinuous and strongly continuous, 
and the function $t\mapsto\mathscr{L}(t)$ is continuous, we immediately obtain 
\begin{eqnarray}
\label{Wliminf}
&{\displaystyle \mathcal{W}(\nabla u_{\infty}(t))
\leq \liminf_{k \to \infty} \mathcal{W}(\nabla u_{k}(t)),\qquad \mathcal{W}(\nabla v)=\lim_{k\to\infty} \mathcal{W}(\nabla v_k), }\\
\label{limL}
&{\displaystyle\langle\mathscr{L}(t),u_{\infty}(t)\rangle=\lim_{k\to\infty}\langle\mathscr{L}_k(t),u_k(t)\rangle, \qquad 
\langle\mathscr{L}(t),v\rangle=\lim_{k\to\infty}\langle\mathscr{L}_k(t),v_k\rangle.}
\end{eqnarray}

So far we have easily obtained that 
\begin{eqnarray}
&&\mathcal{W}(\nabla u_{\infty}(t))-\langle\mathscr{L}(t),u_{\infty}(t)\rangle\leq \nonumber\\ 
&&\qquad\leq \mathcal{W}(\nabla v)-\langle\mathscr{L}(t),v\rangle+
\limsup_{k\to\infty}\|(\phi([v])-\gamma_k(t))^+\|_{1,M}\,,\label{pa0}
\end{eqnarray}
where the last term in right-hand side comes from the equality 
\begin{equation}\label{e+}
(\gamma\vee \phi([v]))-\gamma = (\phi([v])-\gamma)^+,
\end{equation}
which holds for every $\gamma \in L^0(M)^+$.
In order to obtain that the pair $(u_{\infty}(t),\gamma_{\infty}(t))$ satisfies point (a) 
in Definition~\ref{iqe} of quasistatic evolution we want to apply Lemma~\ref{lm_limsup}.
To this aim we need to know that $\phi([u_{\infty}(t)])\in L^1(M)^+$.
By (\ref{bdg}) in Remark~\ref{bounds2} we have that $\|\gamma_k(t)\|_{1,M}$
is bounded uniformly with respect to $k$. 
As $u_k(t)$ belong to $AD(\psi_k(t),\gamma_k(t))$, the sequence $\phi([u_k(t)])$ is bounded in $L^1(M)^+$, and by the
lower semicontinuity assumption $(\varphi_3)$ we obtain that $\phi([u_{\infty}(t)])\in L^1(M)^+$ thanks to the Fatou lemma. 
Then we can apply Lemma \ref{lm_limsup} and we get 
\begin{eqnarray}
&&\mathcal{W}(\nabla u_{\infty}(t))-\langle\mathscr{L}(t),u_{\infty}(t)\rangle\leq\nonumber\\
&& \qquad\leq \mathcal{W}(\nabla v)-\langle\mathscr{L}(t),v\rangle+\|(\phi([v])-\gamma_{\infty}(t))^+\|_{1,M}\,.\label{pa}
\end{eqnarray}
Applying (\ref{e+}) to the last term in the right-hand side of (\ref{pa}) we conclude that
$\mathscr{E}(t)(u_{\infty}(t),\gamma_{\infty}(t))\leq \mathscr{E}(t)(v,\gamma_{\infty}(t)\vee\phi([v]))\leq
\mathscr{E}(t)(v,\gamma)$ for every $t\in[0,T]$ and point (a) of Definition~\ref{iqe} is satisfied.

By the definition of the discrete problems, for every $k$ the function $t\mapsto \gamma_k(t)$ is increasing. 
Passing to the $\sigma^p_\varphi$-limit, the same property holds for $t\mapsto\gamma_{\infty}(t)$, 
so that point (b) of Definition~\ref{iqe} is satisfied.

It remains to prove point (c). For a.e.\ $t$ define
\begin{equation*}
\theta(t):= 
\langle \partial\mathcal{W}(\nabla u_{\infty}(t)),\nabla \dot{\psi}(t)\rangle
-\langle\mathscr{L}(t),\dot{\psi}(t)\rangle-\langle\dot{\mathscr{L}}(t),u_\infty(t)\rangle.
\end{equation*}
Arguing as in the proof of \cite[Theorem~3.15]{DM-F-T-03} we get 
\begin{equation}\label{tinf=t}
\theta_{\infty}(t)=\theta(t),
\end{equation}
for a.e.\ $t\in[0,T]$. This in particular means that the map $t \mapsto \theta(t)$ is measurable.
Since we have proved that for every $t\in[0,T]$ the pair $(u_{\infty}(t),\gamma_{\infty}(t))$ satisfies points 
(a) and (b) of Definition~\ref{iqe}, 
we are in a position to apply Theorem~\ref{th:en>} and get 
\begin{equation*}
\mathscr{E}(t)(u_{\infty}(t),\gamma_{\infty}(t)) - 
\mathscr{E}(0)(u_0,\gamma_0) \geq \int_0^t \theta(s)\,ds.
\end{equation*}
By (\ref{lsc_sp}), (\ref{Wliminf}), and (\ref{limL}) we get
\begin{equation}\label{eq_1}
\mathscr{E}(t)(u_{\infty}(t), \gamma_{\infty}(t))\leq
\liminf_{j\to\infty}\mathscr{E}_{k_j}(t)(u_{k_j}\!(t), \gamma_{k_j}\!(t))\leq
\limsup_{k\to\infty}\mathscr{E}_{k}(t)(u_{k}(t), \gamma_{k}(t)).
\end{equation}
Using Lemma~\ref{energyestim} and taking (\ref{fatoulimsup_tk}) and (\ref{tinf=t}) into account, we obtain 
\begin{equation}\label{eq_2}
\limsup_{k\to\infty}\mathscr{E}_{k}(t)(u_{k}(t), \gamma_{k}(t))\leq
\mathscr{E}(0)(u_0,\gamma_0)+ \int_0^t \theta(s)\,ds.
\end{equation}
By (\ref{eq_1}) and (\ref{eq_2}) we get that
\begin{equation*}
\mathscr{E}(t)(u_{\infty}(t), \gamma_{\infty}(t))\leq
\mathscr{E}(0)(u_0,\gamma_0)+ \int_0^t \theta(s)\,ds
\end{equation*}
holds true for any $t\in[0,T]$, and this concludes the proof.
\end{proof}
In the following theorem we prove that for every $t\in[0,T]$ the energy
for the discrete-time problems converges to the energy for the continuous-time problem. We emphasize that the theorem is true 
for {\em any} irreversible quasistatic evolution $t\mapsto (u(t),\gamma(t))$ corresponding to a given $t\mapsto\gamma(t)$, 
not only for the one obtained as limit of the solutions of the discrete-time problems.
\begin{theorem}\label{thm:d-c}
For every $t\in[0,T]$ let $u_k(t)$ and $\gamma_k(t)$ be defined as in the beginning of the proof of Theorem~\ref{mainthm}.
Assume that $\gamma_k(t)$ $\sigma^p_\varphi$-converges to $\gamma(t)\in L^1(M)^+$ for any $t\in[0,T]$.
Let $t\mapsto (u(t),\gamma(t))$ be an irreversible quasistatic evolution. For a.e.\ $t\in[0,T]$ let $\theta_k(t)$ 
be defined as in (\ref{teta_k}),
and set
\begin{equation*}
\theta(t) := \langle \partial\mathcal{W}(\nabla u(t)),\nabla \dot{\psi}(t)\rangle  -
\langle\mathscr{L}(t),\dot{\psi}(t)\rangle- \langle\dot{\mathscr{L}}(t),u(t)\rangle.
\end{equation*}
Then
\begin{eqnarray}\label{wl}
&{\displaystyle\mathcal{W}(\nabla u(t))-\langle\mathscr{L}(t),u(t)\rangle = \lim_{k\to\infty} 
(\mathcal{W}(\nabla u_k(t))-\langle\mathscr{L}_k(t),u_k(t)\rangle),} \\
&{\displaystyle\|\gamma(t)\|_{1,M} = \lim_{k\to\infty} \|\gamma_k(t)\|_{1,M},}\nonumber
\end{eqnarray}
for every $t\in [0,T]$. Furthermore
\begin{equation*}
\theta_k \to \theta \quad \mbox{in }L^1([0,T]),
\end{equation*}
so that there exists a subsequence of $\theta_k$ which converges to $\theta$ a.e.\ in $[0,T]$.
\end{theorem}
\begin{proof}
For the proof we need to show that
\begin{equation}\label{wlim}
\lim_{j\to \infty}\mathcal{W}(\nabla u_{k_j}(t))= \mathcal{W}(\nabla u_\infty(t)),
\end{equation}
for every $t\in [0,T]$, where $u_{k_j}(t)$ is the subsequence constructed in the proof of Theorem~\ref{mainthm}, 
and $u_\infty(t)$ is its limit.
To this aim, let $v_j:=u_{\infty}(t)+\psi_{k_j}(t)-\psi(t)$. 
By the minimality of the pair $(u_{k_j}(t),\gamma_{k_j}(t))$ we obtain that
$\mathscr{E}_{k_j}(t)(u_{k_j}(t),\gamma_{k_j}(t))\leq \mathscr{E}_{k_j}(t)(v_j,\gamma_{k_j}(t)\vee\phi([u_{\infty}(t)]))$, 
and passing to 
the limit as $j$ goes to infinity, we get by (\ref{eq:egs2}), (\ref{Wliminf}), and (\ref{limL})
\begin{equation}\label{elimsup1}
\begin{split}
&\limsup_{j\to\infty}\big[\mathcal{W}(\nabla u_{k_j}(t))-\langle\mathscr{L}_{k_j}(t), u_{k_j}(t)\rangle\big]\leq\\
& \leq \limsup_{j\to\infty}\big[\mathcal{W}(\nabla v_j)-\langle\mathscr{L}_{k_j}(t), v_j\rangle+
\|(\phi([u_\infty(t)])-\gamma_{k_j}(t))^+\|_{1,M}\big]=\\
&= \mathcal{W}(\nabla u_{\infty}(t))-\langle\mathscr{L}(t), u_{\infty}(t)\rangle+
\limsup_{j\to\infty}\|(\phi([u_\infty(t)])-\gamma_{k_j}(t))^+\|_{1,M}.
\end{split}
\end{equation}
Since $\gamma_{k_j}(t)$ $\sigma^p_\varphi$-converges to $\gamma_{\infty}(t)$, by Lemma $\ref{lm_limsup}$ we have 
\begin{equation}\label{elimsup2}
\limsup_{j \to \infty}\|(\phi([u_{\infty}(t)])-\gamma_{k_j}(t))^+\|_{1,M}\leq 0.
\end{equation}
Taking into account (\ref{elimsup1}) and (\ref{elimsup2}) we get in particular that
\begin{equation}\label{Wlimsup}
\limsup_{j\to\infty}\mathcal{W}(\nabla u_{k_j}(t))\leq \mathcal{W}(\nabla u_\infty(t)).
\end{equation}
This, together with (\ref{Wliminf}), gives~(\ref{wlim}). 

To conclude the proof it is sufficient to follow the arguments of the proof of \cite[Theorem~8.1]{DM-F-T-03}.
\end{proof}

The result can be improved under strict convexity assumption.
\begin{theorem}\label{th:conv}
In addition to the hypotheses of Theorem~\ref{thm:d-c}, assume that $\xi\mapsto W(x,\xi)$ is strictly convex for a.e.\ 
$x\in \Omega\smallsetminus M$ and that $y\mapsto \varphi(x,y)$ is convex for $\mathscr{H}^{n-1}$-a.e.\ $x\in M$. 
Then $u_k(t)\to u(t)$ strongly in $W^{1,p}(\Omega\smallsetminus M;\RR^m)$, for every $t\in [0,T]$.
\end{theorem}
\begin{proof}
We observe that for every $t\in [0,T]$ and $\gamma\in L^1(M)^+$ the functional $v\mapsto \mathscr{E}(t)(v,\gamma)$ 
is strictly convex on the set of functions $v\in W^{1,p}(\Omega \smallsetminus M;\RR^m)$ with 
$v=\psi(t)$ $\mathscr{H}^{n-1}$-a.e.\ on $\partial_0 \Omega$. 
Therefore for every $t$ there exists a unique function 
$u\in AD(\psi(t),\gamma(t))$ such that the pair $(u, \gamma(t))$ is globally stable at time $t$. 
It follows that $u(t)$ coincides with the function $u_\infty(t)$ constructed in the proof of Theorem~\ref{mainthm} 
and that the whole sequence $u_k(t)$ converge to  $u(t)$ weakly in $W^{1,p}(\Omega \smallsetminus M;\RR^m)$. 
Therefore (\ref{wl}) implies that $\mathcal{W}(\nabla u_k(t))\to \mathcal{W}(\nabla u(t))$. 
Using \cite[Theorem 3]{Vis} we deduce that $\nabla u_k(t) \to \nabla u(t)$ in measure. As
\begin{equation*}
|\nabla u_k(t) - \nabla u(t)|^p\leq 2^{p-1}a_0^{-1}[W(\nabla u_k(t))+W(\nabla u(t))]+2^{p-1}a_0^{-1}b_0,
\end{equation*}
the conclusion follows from the generalized dominated convergence theorem. 
\end{proof}  
\end{section}

\begin{section} {Euler conditions}\label{s:euler}

In this section we study the Euler conditions satisfied by globally stable
pairs $(u,\gamma)\in W^{1,p}(\Omega\smallsetminus M;\RR^m)\times L^1(M)^+$. 
Let us fix $t\in [0,T]$ and let $(u,\gamma)\in W^{1,p}(\Omega \smallsetminus M;\RR^m)\times L^1(M)^+$ be globally stable 
at time $t$, and let $v \in W^{1,p}(\Omega \smallsetminus M;\RR^m)$ be such that $v=0$ $\mathscr{H}^{n-1}$-a.e.\ on $\partial_0 \Omega$. 
Hence for every $\varepsilon >0$ the function
$u+\varepsilon v$ belongs to $AD(\psi(t),\gamma\vee\phi([u]+\varepsilon [v]))$, 
and by the global stability of the pair $(u,\gamma)$ at time $t$, we have that
$\mathscr{E}(t)(u,\gamma)\leq 
\mathscr{E}(t)(u+\varepsilon v,\gamma\vee\phi([u]+\varepsilon [v]))$, therefore
\begin{equation}\label{eul_weak}
\liminf_{\varepsilon\to 0^+}
\frac{\mathscr{E}(t)(u + \varepsilon v,\gamma \vee \phi([u] + \varepsilon [v])) -
\mathscr{E}(t)(u,\gamma)}{\varepsilon}\geq 0.
\end{equation}
The weak formulation of the Euler conditions will be obtained from this inequality.
Without loss of generality, we assume that $\mathscr{L}(t)$ is given by (\ref{def_L}),
and we omit the dependence on time. After some standard calculation, one can express (\ref{eul_weak}) in the following form
\begin{equation}\label{eul_w}
\begin{split}
&
\int_{\Omega\smallsetminus M} \big(\partial_\xi W(x,\nabla u) - H\big): \nabla v \, dx - 
\int_{\Omega\smallsetminus M} fv \,dx - \int_{\partial_1 \Omega}gv \, d\mathscr{H}^{n-1} + \\
& - \int_M \big(g^{\oplus}v^{\oplus} +  g^{\ominus}v^{\ominus}\big)\, d\mathscr{H}^{n-1} 
+ \liminf_{\varepsilon\to 0^+}\frac{\|(\phi([u]+\varepsilon [v]) - \gamma)^+\|_{1,M}}{\varepsilon} \, \geq 0,
\end{split}
\end{equation}
for any $v \in W^{1,p}(\Omega \smallsetminus M;\RR^m)$ such that $v=0$ $\mathscr{H}^{n-1}$-a.e.\ on $\partial_0 \Omega$.

To continue our analysis we need now to specify the form of the function $\varphi$. 
More precisely, we consider $\varphi\colon M\times \RR^m\to [0,+\infty]$ defined by
\begin{equation}\label{spec:phi}
\varphi(x,y) := \varphi_0(x) + \tilde{\varphi}(x,y) \quad \text{for }y\neq 0 \quad \text{and}\quad 
\varphi(x,0) := 0\quad\mbox{for all }x\in M, 
\end{equation}
where $\varphi_0\in L^1(M)^+$ and $ \tilde{\varphi}\colon M\times\RR^m\to [0,+\infty]$ is a
Borel function. We assume that for every $x\in M$ the following properties hold:
\begin{itemize}
\item[(1)] $\varphi(x,y)=0$ if and only if $y=0$;
\item[(2)] the function $\tilde{\varphi}(x,\cdot)$ belongs
to the space $C^0(\RR^m)\cap C^1(\RR^m\smallsetminus\{0\})$;
\item[(3)] $\tilde{\varphi}(x,0)=0$;
\item[(4)] there exists an $L^{\infty}$-function $\bar{\varphi}$ such that 
$|\partial_y \tilde{\varphi}(x,y)|\leq \bar{\varphi}(x)$ for any $y\neq 0$, 
where $\partial_y \tilde{\varphi}(x,y)$ denotes the vector of the partial derivatives of $\tilde{\varphi}$ with respect to $y$; 
\item[(5)] the limit
\begin{equation}\label{psit}
\tilde{\psi}(x,y) := \lim_{\varepsilon \to 0^+}\partial_y \tilde{\varphi}(x,\varepsilon y)y  
\end{equation}
exists and is finite for any $y \neq 0$. 
\end{itemize}
\begin{remark}\label{r:Hospital}
By using de l'H${\rm \hat{o}}$pital Theorem, one obtain immediately that 
\begin{equation*}
\tilde{\psi}(x,y) = \lim_{\varepsilon \to 0^+} \frac{\tilde{\varphi}(x,\varepsilon y)}{\varepsilon},
\end{equation*}
for any $x\in M$, $y\neq 0$. It follows from the positiveness of $\tilde{\varphi}$ that $\tilde{\psi}\geq 0$. 
Moreover, we get easily that $\tilde{\psi}$ is positively $1$-homogeneous 
with respect to $y$, i.e., $\tilde{\psi}(x,\lambda y) = \lambda \tilde{\psi}(x,y)$, for every $\lambda>0$.
Furthermore, by (\ref{psit}) and (4), we get also
\begin{equation}\label{psiefibar}
|\tilde{\psi}(x,y)| \leq \bar{\varphi}(x)|y| \quad \mbox{for every $x\in M$ and $y\neq 0$.}
\end{equation}
\end{remark}

The main result of this section is a theorem which makes explicit the Euler conditions obtained from (\ref{eul_w}) 
in the case of the function $\varphi$ specified above.
Before stating the theorem, we establish a general result concerning closed linear subspaces of $L^1_\mu(\Omega)$, for
an arbitrary Radon measure $\mu$ on $\Omega$. We will apply this result to the measure $\mu=\mathscr{H}^{n-1}\LLL M$.

The characteristic function of any set $E$ is denoted by $1_E$, i.e., 
$1_E(x) = 1$ if $x\in E$, $1_E(x) = 0$ otherwise. 
\begin{lemma}\label{lm:linsub}
Let $\mu$ be a Radon measure in $\Omega$ and
let $Y$ be a closed linear subspace of $L^1_\mu(\Omega)$ with the following properties:
\begin{itemize}
\item[(a)] if $u,v\in Y$, then $u\vee v\in Y$;
\item[(b)] if $u\in Y$ and $\omega\in C^\infty_c(\Omega)$, then $\omega u\in Y$.
\end{itemize}
Then there exists a Borel set $E\subset\Omega$ such that $Y=\{{u\in L^1_\mu(\Omega)}: {u=0} \ {\mu\hbox{-a.e.\ on } E}\}$.
\end{lemma}
\begin{proof}
We begin by proving that 
\begin{equation}\label{eq:1}
\mbox{if }u\in L^1_\mu(\Omega) \mbox{ and } |u|\leq |v|  \mbox{ for some }v\in Y, \mbox{ then }u\in Y.
\end{equation}
Indeed in this case there exists $\omega\in L^\infty_\mu(\Omega)$ such that $u=\omega v$ and there 
is a sequence $\omega_k\in C^\infty_c(\Omega)$ such that $\omega_k$ is bounded in $L^\infty_\mu(\Omega)$ 
and $\omega_k \to \omega$ $\mu$-a.e.\ on $\Omega$. 
By (b) we have $\omega_k v\in Y$, and
by the Lebesgue dominated convergence theorem $\omega_k v \to \omega v=u$ in $L^1_\mu(\Omega)$. 
Since $Y$ is closed, we conclude that $u\in Y$.

Now we prove that
\begin{equation}\label{eq:2}
\mbox{if }u\in Y \mbox{ and }t>0, \mbox{ then } u\wedge t\in Y \mbox{ and }(u-t)^+\in Y.
\end{equation}
As $|u\wedge t|\leq |u|$, we have $u\wedge t\in Y$ by (\ref{eq:1}). Since $(u-t)^+=u-u\wedge t$, we obtain that $(u-t)^+\in Y$.

Next we prove that
\begin{equation}\label{eq:3} 
\mbox{if }u\in Y \mbox{ and } t>0, \mbox{ then }1_{\{u>t\}}\in Y, 
\end{equation}
where $\{u>t\}:=\{x\in \Omega:u(x)>t\}$.
By (\ref{eq:2}) we deduce that for every $k>0$ we have $k(u-t)^+\wedge 1\in Y$. As
$ [k(u-t)^+]\wedge 1 \to 1_{\{u>t\}}$ pointwise and $[k(u-t)^+]\wedge 1\leq |u|/t$, 
the convergence takes place in $L^1_\mu(\Omega)$ and we conclude that $1_{\{u>t\}}\in Y$.

Let $(u_k)$ be a sequence dense in $Y$ and let $E$ be the intersection of the sets ${\{u_k=0\}}$. 
It is easy to prove by approximation that $u=0$ $\mu$-a.e.\ on $E$ for every $u\in Y$. 
Conversely, let $u\in L^1_\mu(\Omega)$ with $u=0$ $\mu$-a.e.\ on $E$. For every $k$ let
\begin{equation*}
A_k:=\{u_1\vee u_2 \vee \dots \vee u_k>1/ k\}.
\end{equation*}
By (a) and (\ref{eq:3}) we have $1_{A_k}\in Y$, so that $(k1_{A_k})\wedge u^+$ and $(k1_{A_k})\wedge u^-$ belong to $Y$, 
by (\ref{eq:1}). 
As $(k1_{A_k})\wedge u^+\to u^+$ and $(k1_{A_k})\wedge u^- \to u^-$ in $L^1_\mu(\Omega)$ we conclude that $u\in Y$.
\end{proof}
\begin{lemma}\label{lm:2}
Let $D\subset M$, let $Y_D^m$ be the set of all functions of the form $[v]$, 
with $v\in W^{1,p}(\Omega\smallsetminus M;\RR^m)$
and $[v]=0$ $\mathscr{H}^{n-1}$-a.e.\ on $D$, and let $\overline Y{}_{\!D}^m$ be the closure of $Y_D^m$ 
in $L^1(M\smallsetminus \partial M;\RR^m)$. 
Then there exists a Borel set $\tilde{D}$ (unique up to $\mathscr{H}^{n-1}$-equivalence), containing $D$, 
such that $\overline Y{}_{\!D}^m=\{w\in L^1(M\smallsetminus \partial M;\RR^m): {w=0} \ {\mathscr{H}^{n-1}\hbox{-a.e.\ on }\tilde{D}}\}$.
\end{lemma}
\begin{proof}
Let $Y_D$ be the set corresponding to the case $m=1$.  It is easy to see that $Y_D^m=(Y_D)^m$. 
Therefore it suffices to prove the lemma in the case $m=1$.

The conclusion follows from Lemma~\ref{lm:linsub} applied to $\overline Y{}_{\!D}$. 
It is enough to verify that conditions (a) and (b) are satisfied by 
$Y_D$. Condition (b) is trivial. To prove (a) we consider an open set $U\subset \Omega\smallsetminus M $,  
with $C^1$ boundary and $M\subset \partial U$, such that  $U$ lies on the negative side of $M$. 
Given two functions $u$ and $v \in W^{1,p}(\Omega\smallsetminus M)$ it is easy to check 
that $[u]\vee [v]=[u\vee(v-\tilde{v}+\tilde{u})]$, where $\tilde{u}$ and $\tilde{v}\in W^{1,p}(\Omega)$ 
coincide with $u$ and $v$ on $U$, respectively.
\end{proof}
In the following theorem we will consider a function $u\in W^{1,p}(\Omega \smallsetminus M;\RR^m)$ 
such that the divergence of the matrix field $\partial_{\xi}W(x, \nabla u) - H$ belongs
to $L^q(\Omega \smallsetminus M;\RR^m)$. It turns out that its normal trace $(\partial_{\xi}W(x, \nabla u) - H)\nu$ is defined
as an element of $(W^{1-\frac{1}{p},p}(\partial_1\Omega;\RR^m))'$. Moreover, we have that the normal traces 
$(\partial_{\xi}W(x, \nabla u) - H)^\oplus \nu$ and 
$(\partial_{\xi}W(x, \nabla u) - H)^\ominus \nu$ (defined on the positive and negative side of $M$) are both elements
of the space $(W^{1-\frac{1}{p},p}(M\smallsetminus\partial M;\RR^m))'$.
The duality pairing between $(W^{1-\frac 1 p,p}(M\smallsetminus \partial M;\RR^m))'$ and  
$W^{1-\frac 1 p,p}(M\smallsetminus \partial M;\RR^m)$ will be denoted by $\langle \cdot,\cdot \rangle$. 
\begin{theorem}\label{th:euler}
Let $t\in[0,T]$ and $(u,\gamma)\in W^{1,p}(\Omega \smallsetminus M;\RR^m)\times L^1(M)^+$ 
be globally stable at time $t$. Assume that $\varphi\colon M\times\RR^m\to [0,+\infty]$ is defined as above in (\ref{spec:phi})
and it satisfies (1)--(5). 
Then
\begin{eqnarray}\label{eul_el}
&- {\dv} \big(\partial_{\xi}W(x, \nabla u) - H) = f \quad \text{ on }\Omega \smallsetminus M,\\
\label{eul_neu}
&(\partial_{\xi} W(x,\nabla u)-H)\nu = g \quad \text{ on }\partial_1\Omega,\\
\label{eul_onM}
&(\partial_{\xi} W(x,\nabla u)-H)^{\oplus}\nu + g^{\oplus} = (\partial_{\xi} W(x,\nabla u)-H)^{\ominus}\nu - g^{\ominus}
\quad \text{ on } M\smallsetminus \partial M.
\end{eqnarray}
Let us define 
\begin{eqnarray*}
 A    &:=& \{x\in M : 0<\phi([u])(x)=\gamma(x)\}, \\
 B    &:=& \{x\in M : 0=\phi([u])(x) \mbox{ and } \gamma(x)=\varphi_0(x)\},\\
 D    &:=& \{x\in M : \gamma(x)<\varphi_0(x)\},
\end{eqnarray*}
and let $\tilde{D}$ be the set associated with $D$ by Lemma~\ref{lm:2}.
Then there exists $h\in L^{\infty}({M\smallsetminus \tilde{D}};\RR^m)$ such that 
\begin{equation}\label{hv}
\langle (\partial_{\xi} W(x,\nabla u)-H)^{\oplus}\nu + g^{\oplus},[v]\rangle= 
\int_{M\smallsetminus \tilde{D}} h[v] \,d\mathscr{H}^{n-1},
\end{equation}
for every $v\in W^{1,p}(\Omega \smallsetminus M;\RR^m)$ such that $[v]=0$ $\mathscr{H}^{n-1}$-a.e.\ on $D$.
Moreover 
\begin{itemize}
\item[(a)] for $\mathscr{H}^{n-1}$-a.e.\ $x\in A\smallsetminus \tilde{D}$ the vector $h(x)$
belongs to the segment joining $0$ and $\partial_y\tilde{\varphi}(x,[u](x))$;
\item[(b)] for $\mathscr{H}^{n-1}$-a.e.\ $x\in B\smallsetminus \tilde{D}$ the vector $h(x)$
belongs to the bounded convex set 
$K(x):=\{a\in \RR^m: ay \leq\tilde{\psi}(x,y),\; \forall y\in\RR^m\}$;
\item[(c)] for $\mathscr{H}^{n-1}$-a.e.\ $x\in M\smallsetminus (A\cup B\cup \tilde{D})$ we have $h(x) = 0$.
\end{itemize}
\end{theorem}
\begin{remark}\label{r:D}
It is easy to see that, if $D$ is ($\mathscr{H}^{n-1}$-equivalent to) a closed set, 
then $\tilde{D}=D$ (up to $\mathscr{H}^{n-1}$-equivalence). 
A more difficult proof shows that the same result is true if $D$ is ($\mathscr{H}^{n-1}$-equivalent to) a 
quasi closed set with respect to $(1,p)$-capacity. 

It is clear that, if $\varphi_0=0$, then $\tilde{D}=D=\emptyset$. 
\end{remark}

\begin{remark}\label{r:mecc}
For $\mathscr{H}^{n-1}$-a.e.\ $x\in M$ the vector $h(x)$, obtained in Theorem~\ref{th:euler}, represents 
the cohesive force exerted from the positive lip of the crack on the negative lip. 
The theorem shows the conditions satisfied by the cohesive force on the different regions of $M$ 
determined by the respective relations between $\phi([u])$, $\gamma$ and $\varphi_0$. 
\end{remark}

\begin{proof}[Proof of Theorem~\ref{th:euler}]
Since $\phi([u])\leq \gamma$ $\mathscr{H}^{n-1}$-a.e.\ on $M$, we have
$(\phi([u]) - \gamma)^+=0$ $\mathscr{H}^{n-1}$-a.e.\ on $M$. 
If $[v]=0$ $\mathscr{H}^{n-1}$-a.e.\ on $M$, then the $\liminf$ in (\ref{eul_w}) is actually a limit and it is zero. 
Therefore (\ref{eul_el}), (\ref{eul_neu}), 
and (\ref{eul_onM}) can be obtained from (\ref{eul_w}) by standard argument involving integration by parts and a suitable choice of the 
test function $v\in W^{1,p}(\Omega;\RR^m)$.

To shorten the notation, we set 
$\tilde{h}:= (\partial_{\xi} W(x,\nabla u)-H)^{\oplus}\nu + g^{\oplus}$ on $M\smallsetminus \partial M$. 
As explained before the statement of the theorem, we have 
$\tilde{h}\in (W^{1-\frac 1 p,p}(M\smallsetminus \partial M);\RR^m))'$.
So far, we may rewrite (\ref{eul_w}) as
\begin{equation}\label{euler}
\langle -\tilde{h}, [v]\rangle +
\liminf_{\varepsilon\to 0^+} \frac{\|(\phi([u]+\varepsilon [v])-
\gamma)^+\|_{1,M}}{\varepsilon}\,\geq 0,
\end{equation}
for any $v \in W^{1,p}(\Omega \smallsetminus M;\RR^m)$ such that $v=0$ $\mathscr{H}^{n-1}$-a.e.\ on $\partial_0 \Omega$. 

Let us extend the definition of $\tilde{\psi}$ by setting $\tilde{\psi}(x,0)=0$ for every $x\in M$.
Now we prove that
\begin{equation}\label{lim}
\begin{split}
&\lim_{\varepsilon\to 0^+}\frac{\|(\phi([u]+\varepsilon w) - \gamma)^+\|_{1,M}}{\varepsilon} =\\
& = \int_M \Big((\partial_y \tilde{\varphi}(x,[u])w)^+ \, 1_A + \tilde{\psi}(x,w) \, 1_B\Big) \,d\mathscr{H}^{n-1},
\end{split}
\end{equation}
for every $w\in L^1(M\smallsetminus \partial M;\RR^m)$ with $w=0$ $\mathscr{H}^{n-1}$-a.e.\ on $D$.
To this aim, it is convenient to split the set $M$ into the union of the following 
two disjoint subsets $A':=\{x\in M: [u](x) \neq 0\}$ and $B':=\{x\in M: [u](x) = 0\}$. 

On $A'$, as $\phi([u])\leq \gamma$ $\mathscr{H}^{n-1}$-a.e.\ on $M$, we have that 
\begin{equation*}
\begin{split}
\frac{(\phi([u]+\varepsilon w) - \gamma)^+}{\varepsilon} &\leq \frac{(\phi([u]+\varepsilon w) - \phi([u]))^+}{\varepsilon}
= \frac{(\tilde{\varphi}(x,[u]+\varepsilon w) - \tilde{\varphi}(x,[u]))^+}{\varepsilon} \leq \\
&\leq(\bar{\varphi}(x)w)^+, 
\end{split}
\end{equation*}
$\mathscr{H}^{n-1}$-a.e.\ on $M$, where we used (\ref{spec:phi}), and assumptions (3) and (4). Moreover, we have that 
\begin{equation*}
\frac{(\phi([u]+\varepsilon w)-\gamma)^+}{\varepsilon} \to (\partial_y \tilde{\varphi}(x,[u])w)^+\, 1_A \quad 
\mathscr{H}^{n-1}\mbox{-a.e.\ on } A',
\end{equation*}
because $A= \{ 0 < \phi([u]) = \gamma\}$. 
By the Lebesgue dominated convergence theorem we get
\begin{equation}\label{eq:A'}
\int_{A'} \frac{(\phi([u]+\varepsilon w) - \gamma)^+}{\varepsilon} \,d\mathscr{H}^{n-1}\to 
\int_M (\partial_y \tilde{\varphi}(x,[u])w)^+ \, 1_A \,d\mathscr{H}^{n-1},
\end{equation}
as $\varepsilon \to 0^+$, for every $w\in L^1(M\smallsetminus \partial M;\RR^m)$. 

Let us consider now the integral over $B'$. If $w\in L^1(M\smallsetminus \partial M;\RR^m)$ 
and $w=0$ $\mathscr{H}^{n-1}$-a.e.\ on $D$, we have
\begin{equation*}
\frac{(\phi(\varepsilon w)-\gamma)^+}{\varepsilon}=0 \quad \mathscr{H}^{n-1}\mbox{-a.e. on }D,
\end{equation*}
thus we can focus on the set $B'\smallsetminus D$. 
As $\gamma \geq \varphi_0$ $\mathscr{H}^{n-1}$-a.e.\ on $M\smallsetminus D$, 
for every $w\in L^1(M\smallsetminus \partial M;\RR^m)$ with $w=0$ $\mathscr{H}^{n-1}$-a.e.\ on $D$, we obtain 
\begin{equation*}
\begin{split}
\frac{(\phi(\varepsilon w) - \gamma)^+}{\varepsilon} &\leq \frac{(\phi(\varepsilon w) - \varphi_0)^+}{\varepsilon}
= \frac{\tilde{\varphi}(x,\varepsilon w)}{\varepsilon} \leq \bar{\varphi}(x)|w|
\end{split}
\end{equation*}
$\mathscr{H}^{n-1}$-a.e.\ on $M\smallsetminus D$, where we used (\ref{spec:phi}), and assumptions (3) and (4).
Moreover, by Remark~\ref{r:Hospital} we get that
\begin{equation*}
\frac{(\phi(\varepsilon w)-\gamma)^+}{\varepsilon} \to (\tilde{\psi}(x,w))^+ \, 1_B= \tilde{\psi}(x,w)1_B\quad 
\mathscr{H}^{n-1}\mbox{-a.e.\ on } B',
\end{equation*}
as $\varepsilon \to 0^+$, for every $w\in L^1(M\smallsetminus \partial M;\RR^m)$ with $w=0$ $\mathscr{H}^{n-1}$-a.e.\ on $D$.
We can apply again the Lebesgue dominated convergence theorem and obtain
\begin{equation*}
\int_{B'} \frac{(\phi([u]+\varepsilon w) - \gamma)^+}{\varepsilon} \,d\mathscr{H}^{n-1}\to 
\int_M \tilde{\psi}(x,w)\, 1_B\, d\mathscr{H}^{n-1},
\end{equation*}
as $\varepsilon \to 0^+$, for every $w\in L^1(M\smallsetminus \partial M;\RR^m)$ with $w=0$ $\mathscr{H}^{n-1}$-a.e.\ on $D$.
This concludes the proof of (\ref{lim}). We note that this equality cannot be true if the condition $w=0$ $\mathscr{H}^{n-1}$-a.e.\ 
on $D$ is violated, because in this case 
 \begin{equation*}
\lim_{\varepsilon \to 0^+} \frac{\|(\phi(\varepsilon w)-\gamma)^+\|_{1,M}}{\varepsilon}
= \lim_{\varepsilon \to 0^+} \frac{\|(\varphi_0 + \tilde{\varphi}(\varepsilon w)-\gamma)^+\|_{1,M}}{\varepsilon} = +\infty.
\end{equation*}
Let $Y_D^m$ be the space defined in Lemma~\ref{lm:2}. Notice that $Y_D^m\subset W^{1-\frac 1 p,p}(M\smallsetminus \partial M);\RR^m)$.
By (\ref{euler}) and (\ref{lim}) we have
\begin{equation}\label{eul_ex}
\langle -\tilde{h}, w\rangle 
+ \int_M\big[(\partial_y \tilde{\varphi}(x,[u])w)^+ \, 1_A +
 \tilde{\psi}(x,w)\, 1_B\big]\, d\mathscr{H}^{n-1} \geq 0,
\end{equation}
for any $w\in Y_D^m$.
In order to localize this inequality, we prove first (\ref{hv}).
Due to our assumption (4) and to (\ref{psiefibar}), if we apply (\ref{eul_ex}) to $w$ and $-w$ we deduce that
\begin{equation}\label{abound}
|\langle \tilde{h}, w\rangle| \leq \|\bar{\varphi}\|_\infty 
\|w\|_{1,M\smallsetminus D}\, ,
\end{equation}
for every $w\in Y_D^m$.
It follows that there exists a function
$h\in L^\infty(M\smallsetminus D;\RR^m)$ such that
\begin{equation*}
\langle \tilde{h},w\rangle =\int_{M\smallsetminus D} hw \, d\mathscr{H}^{n-1},
\end{equation*}
for every $w\in Y_D^m$. This implies that (\ref{hv}) is satisfied.
By density from (\ref{eul_ex}) we obtain
\begin{equation}\label{e_x}
\int_{M\smallsetminus D} \Big[ -hw + (\partial_y \tilde{\varphi}(x,[u])w)^+ \, 1_A +
 \tilde{\psi}(x,w) \, 1_B\Big]\, d\mathscr{H}^{n-1} \geq 0,
\end{equation}
for every $w\in \overline{Y}{}_{\!D}^m$. 
Since by Lemma~\ref{lm:2} we have $\overline{Y}{}_{\!D}^m=\{w\in L^1(M\smallsetminus \partial M;\RR^m): {w=0} 
\ {\mathscr{H}^{n-1}\hbox{-a.e.\ on }\tilde{D}}\}$, we conclude that 
\begin{equation}\label{eul_loc}
-h(x) y  + (\partial_y \tilde{\varphi}(x,[u](x))y)^+ \, 1_A(x) + \tilde{\psi}(x,y)\, 1_B(x)\geq 0, 
\end{equation}
for every $y\in \RR^m$ and for $\mathscr{H}^{n-1}$-a.e.\ $x\in M\smallsetminus \tilde{D}$.

In particular, for $\mathscr{H}^{n-1}$-a.e.\ $x\in A\smallsetminus \tilde{D}$ the equality $\partial_y \tilde{\varphi}(x,[u](x))y=0$ 
implies that $h(x)y=0$ (it is enough to use (\ref{eul_loc}) with $y$ and $-y$),
so that for a given $x\in  A\smallsetminus \tilde{D} $ the two vectors $\partial_y \tilde{\varphi}(x,[u](x))$ and $h(x)$ are parallel, 
hence there exists $\lambda(x)$ 
such that 
\begin{equation}\label{eul_M1}
h(x) = \lambda(x) \,\partial_y \tilde{\varphi}(x,[u](x)) \quad 
\mbox{for }\mathscr{H}^{n-1}\text{-a.e.\  }x\in  A\smallsetminus \tilde{D},
\end{equation}
and it is easy to verify that $0 \leq \lambda(x) \leq 1$, by using again (\ref{eul_loc}). In this way we get condition~(a).

On $B\smallsetminus \tilde{D}$, from (\ref{eul_loc}) we obtain
\begin{equation}\label{eul_M2}
- h(x) y + \tilde{\psi}(x,y) \geq 0 \quad
\mbox{for }\mathscr{H}^{n-1}\text{-a.e.\ }x\in B\smallsetminus \tilde{D},
\end{equation}
for every $y\in\RR^m$, which is precisely condition~(b), by the definition of $K$. 
On the remaining part of $M\smallsetminus \tilde{D}$, from (\ref{eul_loc}) we get condition~(c). This concludes the proof.
\end{proof}

\begin{remark}\label{r:onB}
If $\varphi_0(x)>0$ for $\mathscr{H}^{n-1}$-a.e.\ $x\in M$, and $(u,\gamma)= (u(t),\gamma(t))$ 
for an irreversible quasistatic evolution, then (\ref{esssup}) implies that the set $B\smallsetminus \tilde{D}$ 
is nonempty only if there exists $y\in \RR^m\smallsetminus\{0\}$ such that $\tilde{\varphi}(x,y)=0$, for some $x\in M$. 
This happens, for instance, in the Griffith model, where $\varphi$ is given by (\ref{fi_ex}) with $a>0$ and $b=0$.
In this special case, condition (b) becomes $h(x)=0$ $\mathscr{H}^{n-1}$-a.e.\ on $B\smallsetminus \tilde{D}$, because $K(x)=\{0\}$.
\end{remark}

\begin{remark}\label{rem:suff} If for every $x$ the functions 
$\xi\mapsto W(x,\xi)$ and $y\mapsto \varphi(x,y)$ are convex, then
for any $t\in[0,T]$ and $\gamma\in L^1(M)^+$, the functional $u \mapsto \mathscr{E}(t)(u, \gamma\vee \varphi([u]))$ is convex. 
Therefore, it is possible to prove by standard arguments that conditions (a), (b), and (c) of Theorem~\ref{th:euler} 
are equivalent to the inequality
\begin{equation*}
-\int_M hw\, d\mathscr{H}^{n-1} + \lim_{\varepsilon \to 0^+} \frac{\|(\phi([u]+\varepsilon w)-\gamma)^+\|_{1,M}}{\varepsilon}\,\geq 0 ,
\end{equation*}
for every $w\in \overline{Y}{}_{\!D}^m$. Thus, Euler conditions
(\ref{eul_el}), (\ref{eul_neu}), (\ref{eul_onM}), (a), (b), (c) are not only necessary, but also sufficient to global stability. 
\end{remark}

We show now an example of a scalar problem, where the Euler conditions of Theorem~\ref{th:euler} lead to a simplified 
set of boundary conditions.
\begin{example}
Let $m=1$, $p=2$, $W(x,\xi):=\frac 1 2|\xi|^2$, $H(t):=0$, $g^\oplus(t)=g^\ominus(t):=0$, $\phi(y):=|y|$, which correspond to 
the energy functional:
\begin{equation*}
\mathscr{E}(t)(u,\gamma):= \frac 1 2\int_{\Omega\smallsetminus M} |\nabla u|^2\, dx + \int_M \gamma \, d\mathscr{H}^{n-1} - 
\int_{\Omega\smallsetminus M} f(t)u\,dx -\int_{\partial_1 \Omega} g(t)u \, d\mathscr{H}^{n-1}.
\end{equation*}
Let $t\in[0,T]$ and 
$(u,\gamma)\in W^{1,2}(\Omega\smallsetminus M)\times L^1(M)^+$ be globally stable at time $t$.
Then we are in a position to apply Theorem~\ref{th:euler} and the final part of Remark~\ref{r:D}, obtaining 
\begin{equation*}
\begin{cases}
-\Delta u =f(t) & \mbox{on }\Omega\smallsetminus M,\\[3mm] 
u=\psi(t) & \mbox{on } \partial_0 \Omega,\\[3mm] 
{\displaystyle\frac{\partial u}{\partial \nu} =g(t)} & \mbox{on } \partial_1 \Omega,\\[3mm]
{\displaystyle\frac{\partial u}{\partial \nu} =0} & \mbox{on } M\cap \{0\leq|[u]|<\gamma\},\\[3mm]
\Big|{\displaystyle \frac{\partial u}{\partial \nu}}\Big| \leq 1 \mbox{ and } {\displaystyle \frac{\partial u}{\partial \nu}}[u]\geq 
0 &\mbox{on } M\cap \{|[u]|=\gamma\}.\\[3mm]
\end{cases}
\end{equation*}
By Remark~\ref{rem:suff} we have also that if $u$ solves the previous boundary value problem for a given $\gamma$, then the pair 
$(u,\gamma)$ is globally stable at time $t$.
\end{example}

\end{section}

\begin{section}{The case of linear elasticity}\label{linel}

In this section we show that, with some modifications, it is possible to consider also the case where 
the uncracked part of the body is linearly elastic, which is excluded by the first inequality in (\ref{Whp}).

Let $p=2$ and $m=n\geq 1$.  
We assume now that the bulk energy relative to the displacement $u\in W^{1,2}(\Omega\smallsetminus M;\RR^n)$ 
has the form of linear elasticity
\begin{equation*}
\int_{\Omega\smallsetminus M} A(x)Eu{\,:\,}Eu\,dx,
\end{equation*}
where $ Eu:=\frac1 2 (\nabla u+(\nabla u)^T) $ is the symmetric part of the gradient of $u$, 
and $A$ satisfies the following properties:
\begin{itemize}
\item[$(E_1)$]  for every $x\in \Omega$, $A(x)$ is a linear symmetric operator from the space $\MM^{n\times n}_{sym}$ 
of symmetric $n\times n$ matrices into itself, and the map $x\mapsto A(x)$ is measurable;
\item[$(E_2)$] there are two positive constants $c_0$ and $c_1 $ such that  
\begin{equation}\label{Ehp}
c_0 \left |\xi\right|^{2}\leq A(x)\xi : \xi \leq c_1 \left |\xi\right|^{2} 
\end{equation} 
for every $x \in \Omega\smallsetminus M$ and $\xi \in \MM^{n \times n}_{sym}$.
\end{itemize}

For the sake of simplicity in the notation we introduce the $C^1$ map
$\mathcal{Q}\colon L^2(\Omega\smallsetminus M; \MM^{n\times n}_{sym})\to \RR$ defined by
\begin{equation*}
\mathcal{Q}(\Psi):= \int_{\Omega\smallsetminus M} A(x)\Psi:\Psi\, dx
\end{equation*}
for every $\Psi \in L^2(\Omega\smallsetminus M; \MM^{n\times n}_{sym})$, whose 
differential 
$\partial\mathcal{Q}\colon L^2(\Omega\smallsetminus M; \MM^{n\times n}_{sym})\to L^2(\Omega\smallsetminus M; \MM^{n\times n}_{sym})$
is given by
\begin{equation*}
\langle \partial\mathcal{Q}(\Psi),\Phi\rangle = 
2\int_{\Omega\smallsetminus M} A(x)\Psi{\,:\,}\Phi\;dx,
\end{equation*}
for every $\Phi$, $\Psi\in L^2(\Omega\smallsetminus M;\MM^{n\times n}_{sym})$, where 
$\langle \cdot,\cdot\rangle$ 
denotes now the scalar product in the space $L^2(\Omega\smallsetminus M;\MM^{n\times n}_{sym})$.

For every $t\in[0,T]$ the total energy of an admissible configuration 
$(u,\gamma)\in W^{1,2}(\Omega\smallsetminus M, \RR^n)\times L^1(M)^+$ at time $t$ is now defined as
\begin{equation*}
\mathscr{E}(t)(u,\gamma):= \mathcal{Q}(Eu) - \langle\mathscr{L}(t),u\rangle + \|\gamma\|_{1,M}.
\end{equation*}
Once we have the energy functional, we introduce the notion of global stability as in Definition~\ref{gs}.

Since the $(n-1)$-dimension of $\partial_0\Omega$
is positive, Korn inequality holds (see, e.g., \cite{Cia-97}, \cite{D-L-76}): 
there exists a constant $C=C(\Omega,\partial_0\Omega)$ such that
\begin{equation*}
\|\nabla u\|_2\leq C\|Eu\|_2 \quad \text{for all $u\in W^{1,2}(\Omega;\RR^n)$ such that $u=0$ on $\partial_0\Omega$.}
\end{equation*}
As an immediate consequence, we get the following Korn-type inequality:
\begin{equation}\label{Eu}
\|\nabla u\|_2\leq C\|Eu\|_2 + (C+1)\|\nabla \psi\|_2\, 
\end{equation}
for every $u\in W^{1,2}(\Omega\smallsetminus M;\RR^n)$, and $\psi\in W^{1,2}(\Omega;\RR^n)$
such that
$u=\psi$ on $\partial_0\Omega$. 

Thanks to (\ref{Eu}), we still have an a priori bound for the displacement $u$ as in Remark~\ref{r:apriori}. 

The definition of irreversible quasistatic evolution of minimum energy configurations is now given replacing 
$\langle \partial\mathcal{W}( \nabla u(t)), \nabla \dot{\psi}(t)\rangle$ by 
$\langle \partial\mathcal{Q}( Eu(t)), E\dot{\psi}(t)\rangle$ in Definition~\ref{iqe}.

Thanks to the Korn-type inequality (\ref{Eu}), Theorems \ref{tm:essup}, \ref{mainthm}, \ref{thm:d-c}, and \ref{th:conv} 
(and Remark~\ref{r:c-c'}) continue to hold, with essentially the same proofs, if we replace $\mathcal{W}( \nabla u(t))$ and 
$\langle \partial\mathcal{W}( \nabla u(t)), \nabla \dot{\psi}(t)\rangle$  by $\mathcal{Q}( Eu(t))$ and 
$\langle \partial\mathcal{Q}( Eu(t)), E\dot{\psi}(t)\rangle$, respectively, and a similar substitution is done for $u_k(t)$.

\end{section}

\end{document}